\documentclass{cmslatex}
          
\usepackage[paperwidth=7in, paperheight=10in, margin=.8in]{geometry}

\usepackage{amsmath,amsfonts,amssymb}
\usepackage{enumitem}
\usepackage{hyperref}
\usepackage{lineno}
\usepackage{color}
\usepackage{graphicx,epsfig}
\usepackage[percent]{overpic}
\usepackage{tcolorbox}
\usepackage{tikz-cd}

\usepackage{subfigure}

\newcommand{\R}{\mathbb{R}}
\newcommand{\alphastar}{\alpha^*}
\newcommand{\rtheta}{\rho^\theta}
\newcommand{\Vi}{V_\textsc{int}}
\renewcommand{\r}{\rho}
\newcommand{\Dx}{\Delta x^1}
\newcommand{\Dy}{\Delta x^2}
\newcommand{\Dt}{\Delta t}

\newcommand{\REV}[1]{\textcolor{black}{#1}}

\renewcommand{\theequation}{\arabic{section}.\arabic{equation}}

\newtheorem{thm}{Theorem}[section]

\newtheorem{rem}[thm]{Remark}

\begin{document}
\title{A Generalized Mean-Field Game Model 	for \\ the  
	Dynamics of Pedestrians with \\ Limited Predictive Abilities
}


\author{
Emiliano Cristiani\thanks{
	Istituto per le Applicazioni del Calcolo ``M.\ Picone'', 
	Consiglio Nazionale delle Ricerche,
	Via dei Taurini, 19 --
	00185 Rome, Italy, 
	e.cristiani@iac.cnr.it (corresponding author)}
\and Arianna De Santo\thanks{
	Sapienza -- Universit\`a di Roma, Rome, Italy}
\and Marta Menci\thanks{
	Istituto per le Applicazioni del Calcolo ``M.\ Picone'', 
	Consiglio Nazionale delle Ricerche, Rome, Italy, 
	m.menci@iac.cnr.it
}
}




\pagestyle{myheadings} 

\markboth{A Generalized Mean-Field Game Model with Limited Predictive Abilities}{E. Cristiani, A. De Santo, M. Menci} 

\maketitle

\begin{abstract}
This paper investigates the model for pedestrian flow firstly proposed in [Cristiani et al., DOI:10.1137/140962413].
The model assumes that each individual in the crowd moves in a known domain, aiming at minimizing a given cost functional. Both the pedestrian dynamics and the cost functional itself depend on the position of the whole crowd. 
In addition, pedestrians are assumed to have predictive abilities, but \emph{limited in time}, extending only up to $\theta$ time units into the future, where $\theta\in[0,\infty)$ is a model parameter. 
1) For $\theta=0$ (no predictive abilities), we recover the modeling assumptions of the Hughes's model, where people take decisions on the basis of the current position of the crowd only. 
2) For $\theta\to\infty$, instead, we recover the standard mean field game (MFG) setting, where people are able to forecast the behavior of the others at any future time and take decisions on the basis of the current and future position of the whole crowd. 
3) For very short values of $\theta$ (typically coinciding with a single time step in a discrete-in-time setting), we recover instead the MFG setting joined to the \REV{instantaneous} model predictive control technique.
4) For intermediate values of $\theta$ we obtain something different: as in the Hughes's model, the numerical procedure to solve the problem requires to run an off-line procedure at any fixed time $t$, which returns the current optimal velocity field at time $t$ by solving an associated backward-in-time Hamilton--Jacobi--Bellman equation; but, differently from the Hughes's model, here the procedure involves a prediction of the crowd behavior in the sliding time window $[t,t+\theta)$, therefore the optimal velocity field is given by the solution of a forward-backward system which joins a Fokker--Planck equation with a Hamilton--Jacobi--Bellman equation as in the MFG approach.
The fact that a different forward-backward system must be solved at any time $t$ arises new interesting theoretical questions. Numerical tests will give some clues about the well-posedness of the problem.
\end{abstract}
          



\section{Introduction.}\label{sec:intro}	
This paper deals with the modeling of pedestrian dynamics in the framework of mean field games, with a particular focus on numerical tests which, hopefully, will give useful insights for future theoretical investigations.

Pedestrian modeling has a long-standing tradition, starting from the pioneering papers by Hirai and Tarui \cite{hirai1975}, Okazaki \cite{okazaki1979TAIJa}, and Henderson \cite{henderson1974} in the '70s.
Since then, all types of models were proposed, spanning from microscale to macroscale, including multiscale ones, both differential and nondifferential (e.g., cellular automata). A number of review papers and books are now available \cite{aghamohammadi2020, 
	bellomo2011,
	cristiani2014book, 
	dong2020, 
	duives2013,
	eftimie2018,
	haghani2020, 
	martinez2017,
    papadimitriou2009,
    rosini2013book}, 
we refer the interested reader to these references for an introduction to the field. 
It is also useful to mention that models for pedestrians often stem from those developed in the context of vehicular traffic \cite{helbing2001, rosini2013book}.
Moreover, there is a strict connection between pedestrian modeling and control theory, see, e.g., \cite{albi2020, cristiani2014book} and reference therein.

Beside the scale of observation, a good criterion to classify pedestrian flow models is on the basis of their \textit{rationality}, \textit{predictive ability} and \textit{knowledge of the environment}.
\begin{itemize}[leftmargin=*]
\item At the top of the list we find the models based on the \textit{mean field game} (MFG) theory \cite{arjmand2021a, arjmand2021b, burger2014, dogbe2010, gueant2011, hoogendoorn2004, lachapelle2011}.  
	In this case it is assumed that each pedestrian perfectly knows the environment it moves in, and is able to forecast the movement of the group mates \textit{at any future time}. 
	Moreover, each pedestrian wants to follow a trajectory which minimizes some cost criterion (e.g., it wants to reach a given target in minimum time). It is also assumed that the dynamics of each pedestrian and/or the cost of each trajectory is influenced by the distribution of the whole crowd, therefore pedestrians moving in the same environment are in competition with each other: each pedestrian cannot select its optimal path until all the others have chosen their path, i.e.\ until the crowd as a whole has defined its dynamics. 
	This ``game" hopefully leads to a Nash equilibrium, i.e.\ a set of paths, one for each person, such that no player can lower its own expected cost by changing its strategy while the other players keep theirs unchanged.
	
	The MFG approach is based on the coupling of a forward-in-time Fokker--Planck (FP) equation, which describes the evolution of the density of the crowd, and a backward-in-time Hamilton--Jacobi--Bellman (HJB) equation, which returns the optimal velocity field to be used to drive the crowd dynamics at any time.
	
	Using the theory of MFG for describing pedestrian dynamics was criticized by some authors (see, e.g., \cite{burger2014}) since it seems that assumptions are not totally realistic. In fact, it is difficult that all pedestrians constituting a crowd have such a predictive abilities, unless we consider a crowd moving in a very well known area in a recurring way like, e.g., daily commuters in a train station.
	\item In order to create more usable models, crowd predictive abilities were reduced to a certain extent. Some authors proposed to include a \emph{discount factor} in the cost functional in such a way that long-term forecast of the crowd distribution has a low impact in the instant decisions, see e.g.\ \cite{bardi2020, bertucci2019, degond2017}. 
	These models, often referred to as ``myopic" or ``partially blind", try to describe the fact that people are able to forecast the behavior of themselves and the others only to a certain extent. However, it is important to note that the discount factor only lowers the weight of future decisions without completely annihilating it.	
	\item If we assume instead that pedestrians are not able to forecast the movements of the crowd we fall in the Hughes-type models. 
	Hughes's model \cite{hughes2002} (see also \cite{huang2009} for its interpretation in the framework of optimal control theory) assumes that pedestrians perfectly know the whole environment and are able to compute their optimal path taking into account the distribution of the whole crowd at the current instant only. 
	They cannot predict the positions of the others, nor, a fortiori, take decisions based on that prediction.
	
	The Hughes's model is based on two coupled PDEs: a forward-in-time conservation law for the evolution of the density of the crowd, and a stationary HJB (in particular, the Eikonal equation), which returns -- at any fixed time -- the optimal velocity field to be used to drive the crowd dynamics at that time.
	
	Hughes's model was deeply investigated from the modeling \cite{burger2014, jiang2009, twarogowska2014}, numerical \cite{carlini2017}, and theoretical point of view \cite{amadori2012, camilli2017, difrancesco2011, eftimie2018, elkhatib2013, goatin2013}.  
    To further reduce the capabilities of the pedestrians, one can assume that they have no clear view of the whole environment. This is done introducing a \emph{cone of vision} so that decisions are made taking into account only the nearest part of the crowd \cite{carrillo2016}.	
	\item Finally, if we assume that pedestrians have a target but they are not able to compute the optimal path to reach it given the crowd distribution, we fall in the framework of basic models. These models are typically constituted by two ingredients: 1) a given (precomputed) velocity/force field, which steers all single pedestrians to the target assuming that no other is present in the area. 2) An interaction velocity/force field which describes the interactions with group mates.
	The interaction field is typically constituted by a repulsive force exerted by the crowd on each pedestrian. It can be \emph{local}, if only the pointwise density of the crowd is used to evaluate the interaction with group mates, or \emph{nonlocal}, if each pedestrian evaluates the density of the crowd in a surrounding area. The sum of the two fields gives the final velocity/force which steer the pedestrian motion, see, e.g., \cite{burstedde2001} for a cellular automata, \cite{coscia2008} for a first-order macroscopic model, \cite{bellomo2011} for a second-order macroscopic model, \cite{colombo2012, cristiani2011} for nonlocal models, \cite{helbing1991, helbing2000, helbing1995} for the classical agent-based social force model.
	These models can show some artifacts. In fact, if some obstacles are present, pedestrians can be pushed by the crowd \emph{inside} the obstacles. A suitable treatment of the obstacles (including walls) must be then adopted, see, e.g., \cite{cristiani2017} for some hints. 
\end{itemize}

\medskip

In this paper we investigate the model for pedestrian flow firstly proposed in \cite{cristiani2015SIAP}. As in the MFG approach, the model assumes that each individual in the crowd knows the domain it is moving in and aims at minimizing a given cost functional. 
Both the pedestrian dynamics and the cost functional itself depend on the position of the whole crowd. 
We consider both the finite-horizon and minimum-time problems. 
The main novelty is that \textit{pedestrians do have predictive abilities, but they are limited in time}, extending only up to $\theta$ time units into the future, where $\theta\in[0,\infty)$ is a model parameter. 
To our knowledge, this is the first paper which considers a customizable time-window for prediction. 
Numerical tests will show the impact of such parameter on crowd dynamics.

Note that for $\theta=0$ (no predictive abilities) we recover the modeling assumptions of the Hughes's model, while for $\theta\to\infty$ we recover the pure MFG setting. 
For very short values of $\theta$ (typically coinciding with a single time step in a discrete-in-time setting), instead, we recover the MFG setting joined with the \REV{instantaneous} model predictive control technique \cite{degond2017}. 

For intermediate values of $\theta$, we get a rather complex behavior of the crowd, which reflects on the numerical algorithm used to solve the equations: as in the Hughes's model, the numerical procedure requires to run an offline procedure at any fixed time $t$, which returns the current optimal velocity field at time $t$ by solving an associated backward-in-time HJB equation; but, differently from the Hughes's model, here the procedure involves a prediction of the crowd behavior in the sliding time window $[t,t+\theta)$, therefore the optimal velocity field is given by the solution of a forward-backward system which joins a forward FP equation with a backward HJB equation as in the MFG approach.
The fact that at any time $t$ a different forward-backward system must be solved arises new interesting theoretical questions. Numerical tests will give some clues about the well-posedness of the problem.

\medskip

The paper is organized as follows:
in Section \ref{sec:model} we present the model in the framework of both finite-horizon and minimum-time problems. In Section \ref{sec:numerics} we give some details about the numerical approximation and the algorithm to compute the solution of the considered equations.
In Section \ref{sec:tests} we present several numerical tests, also discussing some numerical evidences about the well-posedness of the problems. 


\section{The model.}\label{sec:model}
Let us consider a mass of numerous and indistinguishable people moving in a two-dimensional domain, and let us denote by $\r(x,t)$ their density function at point $x$ and time $t$.

\medskip

\subsection{Finite-horizon problem.}
In order to set up the finite-horizon control problem, let us  fix a time horizon $T>0$ and focus on a \textit{single test pedestrian}, denoting its position by $y(t)\in\R^2$.
The pedestrian moves in the space interacting with the surrounding crowd. We describe its dynamics by the following stochastic ordinary differential equation
\begin{equation}\label{dynamics_single_pedestrian}
\left\{\begin{array}{ll}
dy(s)=V\big(y,s;\alpha(y,s),\r(\cdot,s)\big)\ dt+\sqrt{2\sigma}\ dW_s, \qquad t<s\leq T \\ [2mm]
y(t)=x
\end{array}
\right.
\end{equation}
where $t\geq 0$ is a generic initial time, $x\in\R^2$ is a generic initial position, $\sigma\geq 0$ is the diffusion parameter, and $W$ is the Wiener process. The velocity term $V\in\R^2$ is given by the sum of two terms:
$$
V\big(y,s;\alpha(y,s),\r(\cdot,s)\big):=\alpha(y,s)+\Vi\big(y,s;\alpha(y,s),\r(\cdot,s)\big).
$$
The control $\alpha$ is the part of the dynamics which the pedestrian can govern, and can be freely chosen in a set of \textit{admissible controls} $\mathcal A$
$$
\mathcal A:=\{f:\R^2\times[0,T]\to A\ :\ A\subset \R^2\ \text{bounded}\}.
$$ 
The interaction velocity $\Vi$, instead,  takes into account the interactions with the crowd and depends both on $\alpha(y,s)$ and the whole distribution $\r$ at time $s$, possibly in a nonlocal way;
we will denote the solution to \eqref{dynamics_single_pedestrian} by $\bar y=\bar y(s;x,t,\alpha,\r)$, or simply $\bar y(s)$ for brevity.

Let us now define the cost functional $J:\R^2\times[0,T]\times\mathcal A\times L^\infty(\R^2\times[0,T])\to\R^+$ as
\begin{equation}\label{def:J}
J(x,t;\alpha,\r):=\int_t^T \ell \big(\bar y(s),s,\alpha(\bar y(s),s),\r(\cdot,s) \big) \ ds + g(\bar y(T)),
\end{equation}
where $\ell:\R^2\times[0,T]\times A\times L^\infty(\R^2)\to\R^+$ is a given \textit{running cost} and $g:\R^2\to\R^+$ a given \textit{terminal cost}.

The optimal control problem of our interest consists in finding the optimal choice $\alphastar\in\mathcal A$ of the control such that $J$ is minimal, i.e.
\begin{equation}\label{alphastar}
J(x,t;\alphastar,\r)=\min_{\alpha\in\mathcal A} J(x,t;\alpha,\r).
\end{equation}

It is well known \cite[Rem.\ 3.10]{bardibook} that $\alphastar$ can be found by first solving for $\phi$ a backward-in-time HJB equation of the form
\begin{equation}\label{HJB-orizzontefinito}
\hspace{-0.4cm}\left\{\begin{array}{l}
-\partial_t\phi(x,t;\r) + \max\limits_{a\in A}\left\{-V(x,t;a,\r)\cdot\nabla\phi(x,t;\r)-\ell(x,t,a,\r)\right\}=\\ [2mm] \hspace{6cm} 
=\sigma\Delta\phi(x,t;\r),\quad x\in\R^2, \quad t\in (0,T) \\ [0mm]
\phi(x,T)=g(x),\quad x\in\R^2
\end{array}\right.
\end{equation}
where $\nabla$ denotes the gradient w.r.t.\ $x$, and then taking 
\begin{equation}\label{synthesis}
\alphastar(x,t;\r)\in\arg\max\limits_{a\in A}\left\{-V(x,t;a,\r)\cdot\nabla\phi(x,t;\r)-\ell(x,t,a,\r)\right\},\quad x\in\R^2, \quad t\in (0,T).
\end{equation}
Note that more than one optimal control could exist. 
Plugging $\alpha^*$ into \eqref{dynamics_single_pedestrian} and solving the equation, we finally get an \textit{optimal trajectory} $\bar y^*$ for the test pedestrian. 

\medskip

In the context of pedestrian dynamics, it is usual to choose the interaction velocity as a repulsion (social) force which acts against the others. 
In this way pedestrians tend to avoid most crowded regions, since they are uncomfortable or simply because they slow down the motion. 
Similarly to \cite{cristiani2015SIAP}, we define
\begin{equation}\label{def:Vint}
\Vi\big(y,s;\alpha(y,s),\r(\cdot,s)\big):=-\int_{\mathcal S(y;\alpha)}C_\textsc{rep}\frac{\zeta-y}{|\zeta-y|^2} \r(\zeta,s)d\zeta
\end{equation} 
where $C_\textsc{rep}\geq 0$ is a model parameter and
\begin{equation}\label{def:sensoryregion}
\mathcal S(y;\alpha):=\left\{\zeta\in\R^2 \ : \ R_0\leq|\zeta-y|\leq R \text{ and } (\zeta-y)\cdot\alpha>0 \right\}
\end{equation}
is a \textit{sensory region} which defines the zone where the crowd has an influence on the pedestrian dynamics. 
The parameters $R_0, R>0$ rule the size of the sensory region, while, assuming that pedestrians walk pointing their head and gaze in the direction $\alpha$, the condition $(\zeta-y)\cdot\alpha>0$ translates the fact that pedestrians react to what it is in front of them only.

Overall, the dynamics introduced above lead to the following behavior of the pedestrian: it points towards the direction $\alpha$ but it is repulsed by people it has in front through $\Vi$. 
It tries to place itself in the regions of the space where $\ell$ is low and, at final time $T$, it tries to find itself where $g$ is low.

\medskip

\paragraph{Mean-Field Game.}
If we assume that \textit{all} the pedestrians want to behave optimally minimizing the same functional cost $J$, we enter the field of MFGs. 
In fact, the optimal trajectory $\bar y^*$ for a single pedestrian can be computed only if the entire distribution $\r$ of the crowd is known at every time $t\in[0,T]$; but -- in turn -- the distribution of the crowd can be computed only if the motion of each pedestrian constituting the crowd itself is known. This leads to a game where each pedestrian is in competition with each other and need to guess the others' behavior in order to find its own optimal strategy. If all people succeed, i.e.\ everyone behaves optimally assuming the others do not change their strategy unilaterally, the system has reached a \textit{Nash equilibrium}, which is the kind of equilibrium we are interested in. 

Following the fundamentals of the MFG theory, we know that, once the equilibrium is reached, the evolution of the density function $\r$ is described by the following Fokker-Planck equation
\begin{equation}\label{FP}
\left\{\begin{array}{ll}
\partial_s \r(y,s)  + \text{div} \Big(\rho V\big(y,s;\alpha^*(y,s;\r),\r(\cdot,s)\big)\Big) = \sigma\Delta\r(y,s),\qquad y\in\R^2, \quad s\in(0,T] \\ [2mm]
\r(y,0)=\r_0(y),\qquad y\in\R^2
\end{array}\right.
\end{equation}
where $\r_0$ is the initial spatial distribution of the crowd and $\alpha^*$ is found as in \eqref{synthesis} using the solution to \eqref{FP} as input. 
It is clear that \eqref{HJB-orizzontefinito} and \eqref{FP} are coupled together and must be solved as a one in the entire time window $[0,T]$.

\medskip

\paragraph{Limited prediction ability.}
We are now ready to introduce the model we want to investigate. 
The main novelty is that the optimal control problem for the test pedestrian \eqref{dynamics_single_pedestrian}-\eqref{alphastar} is solved assuming that, at any time $s$, the density function $\r$ is known only until time $s+\theta$.
This means that pedestrians forecast the evolution of the crowd only for a time $\theta$ in the future and compute their optimal strategy on the basis of that prediction.
Since all the pedestrians are able to do that, we face again a MFG, but restricted in the time window $[s,s+\theta]$. 
Moreover, since the time window moves forward as time $s$ increases, the game \textit{changes continuously in time} and a new Nash equilibrium must be found \emph{every time $s$}.

Note that the presence of the terminal condition $g$ at time $T$ force us to solve the backward-in-time HJB equation in the whole time horizon $[0,T]$. 
To do that, we need to have $\r$ defined in the whole time horizon too. This can be achieved by prolonging $\r$ after time $s+\theta$ in some way, for example freezing the prediction at time $s+\theta$ as we do in the following.\footnote{Another trivial way to prolong the density in time could be to assume that the density is null from time $s+\theta$ to time $T$. Our choice to freeze the density at $s+\theta$ is inspired by the Hughes's model, in which the optimal control problem is solved assuming the current density persists until final time.}

\newpage 
The model reads as
\begin{tcolorbox}[colback=blue!5!white]
\begin{equation}\label{model:main_forward}
\rotatebox[origin=c]{90}{\tiny \text{\textcolor{red}{forward}}}
\left\{\begin{array}{l}
\partial_s \r(y,s)  + \nabla\cdot \big(\r V(y,s;\alphastar,\r)\big) = \sigma\Delta\r(y,s),\quad y\in\R^2, \quad s\in(0,T) \\ [2mm]
\r(y,0)=\r_0(y),\quad y\in\R^2
\end{array}\right. 
\end{equation}
where, for any fixed $s$, $\alphastar(y,s)$ is computed by solving the following forward-backward system:
\begin{equation}\label{model:mfg}
\hspace{-0.3cm}
\rotatebox[origin=c]{90}{\tiny \text{\textcolor{red}{MFG}}}
\left\{
\begin{array}{l}
\rotatebox[origin=c]{90}{\tiny \text{\textcolor{red}{forward}}}
\left\{\begin{array}{l}
\partial_\tau \tilde \r(y,\tau)  + \nabla\cdot \big(\tilde \r V(y,\tau;\alphastar,\tilde \r)\big) = \\ [1mm]
\hspace{3.5cm}=\sigma\Delta\tilde \r(y,\tau),\quad y\in\R^2, \quad \tau\in\big(s,\min\{s+\theta,T\}\big) 
\\ [2mm]
\tilde \r(y,s)=\r(y,s),\quad y\in\R^2
\end{array}\right. 
\\ [8mm]
\rtheta(y,t):=
\left\{
\begin{array}{lll}
\r(y,t), & t\leq s & \text{\small[already acquired]}\\
\tilde \r(y,t), & s< t \leq \min\{s+\theta,T\} & \text{\small[prediction]}\\
\tilde \r(y,s+\theta), &  s+\theta<t\leq \min\{s+\theta,T\} & \text{\small [freezing]}
\end{array}
\right.
\\ [8mm]
\rotatebox[origin=c]{90}{\tiny \text{\textcolor{red}{backward}}}
\left\{\begin{array}{l}
-\partial_t\phi(x,t;\rtheta) + \max\limits_{a\in A}\left\{-V(x,t;a,\rtheta)\cdot\nabla\phi(x,t;\rtheta)-\ell(x,t,a,\rtheta)\right\}=\\ [2mm] \hspace{4.8cm} 
=\sigma\Delta\phi(x,t;\rtheta),\quad x\in\R^2, \quad t\in (0,T) \\ [0mm]
\phi(x,T)=g(x),\quad x\in\R^2
\end{array}\right.
\\ [8mm]
\alphastar(x,t)\in\arg\max\limits_{a\in A} \left\{-V(x,t;a,\rtheta)\cdot\nabla\phi(x,t;\rtheta)-\ell(x,t,a,\rtheta)\right\},  \\ \hspace{8cm} 
x\in\R^2,\quad t\in[0,T].
\end{array}
\right.
\end{equation}
\end{tcolorbox}
Equation \eqref{model:main_forward} is the main dynamics from initial to final time. 
The density $\r$ solution to this equation is the function we will show in the numerical tests.
The system \eqref{model:mfg}, instead, is the ancillary MFG to be solved at any time $s$. The forward part predicts the dynamics between $s$ and $s+\theta$, then $\tilde \r$ is used to reconstruct the density function $\rtheta$ in the whole time interval $[0,T]$. 
Finally $\rtheta$ is used to computed the optimal control $\alphastar$ through the HJB equation. 
Note that, once $\alphastar(\cdot,\cdot)$ is found, only $\alphastar(\cdot,s)$ is actually used to move forward in \eqref{model:main_forward}. (and yes, it is a pity to throw away all those calculations.)

\medskip

\begin{rem}\label{rem:theroeticalinsights}
The question arises if the model \eqref{model:main_forward}-\eqref{model:mfg} is covered by the classical theory of MFGs, so that well-posedness can be proven using established results \cite{lasry2007}. 
The answer seems to be negative due to the fact that a different MFG has to be solved each time $s$. 
Even if the regularity of $\alphastar$ can be proven at any \textit{fixed} time $s$, it is more difficult to get regularity results with respect to time. 
\end{rem}

\subsection{Minimum-time problem.}
In the context of pedestrian dynamics it is common to consider the minimum-time-to-target problem rather than finite-horizon problem. 
This is due to the fact that in typical evacuation problems people are confined in a bounded domain and they want to reach the exit as far as possible. If multiple exits are present and/or the density is large enough to create congestion, the optimal evacuation strategy can be hard to devise.

To this end, let us denote by $\Omega\subset\R^2$ the bounded domain the pedestrians move in and want to leave. Let us also denote by $\mathcal T\subseteq\partial\Omega$ the \emph{target} (exits) on the boundary of the domain to be reached in minimal time. 
We drop the functions $\ell$ and $g$ and we rewrite the functional cost as the first time the trajectory hits the target
\begin{equation}\label{def:Jmtp}
J(x,t;\alpha,\r):=\min\{s:\bar y(s;x,t,\alpha,\r)\in\mathcal T\}
\end{equation}
(could be $+\infty$ if the trajectory never hits the target).
Finally, we consider a fully deterministic dynamics by setting $\sigma=0$ in \eqref{dynamics_single_pedestrian}.
In the end, the two FP equations in the model \eqref{model:main_forward}-\eqref{model:mfg} remain unchanged (apart from the diffusion term which disappears), while the HJB equation assumes the form
\begin{equation}\label{HJB-mtp}
\hspace{-0.4cm}\left\{\begin{array}{l}
-\partial_t\phi(x,t;\rtheta) + \max\limits_{a\in A}\left\{-V(x,t;a,\rtheta)\cdot\nabla\phi(x,t;\rtheta)\right\}-1=0,\quad x\in\Omega\backslash\mathcal T, \quad t>0 \\ [2mm]
\phi(x,t)=0,\quad x\in\mathcal T, \quad t\geq 0 \\ [2mm]
\phi(x,t)=+\infty,\quad x\in\partial\Omega\backslash\mathcal T, \quad t\geq 0
\end{array}\right.
\end{equation}
see \cite[Sect.\ 4.4]{cristiani2014book} for a more detailed derivation.

In this context, the value function $\phi(x,t)$ indicates the minimal time to reach the target starting from $x$ at time $t$. The boundary conditions on $\partial\Omega$ is added so that it is never convenient to reach the boundary of the domain (unless an exit is there). 
In this way, the optimal velocity field will tend to keep the pedestrians inside the domain. 
However, the presence of the interaction velocity $\Vi$, which is out of control of the pedestrians, could lead to a velocity field which actually points outside the domain. Therefore, an impermeable boundary condition is needed in the FP equation too. 
To enforce this constraint we modify the velocity field $V$, zerofying the component of the vector $V$ which points outside the domain.

\section{Numerical approximation and computation.}\label{sec:numerics}
The numerical approximation of \eqref{model:main_forward}-\eqref{model:mfg} is challenging. 
In order to avoid an uncontrolled propagation of the numerical error,  we employ simple first-order schemes for both FP and HJB equations. 

In order to actually implement a numerical scheme on a finite grid, a final time $T$ is necessarily introduced also for the minimum-time problem. In this case, we manually check that the entire crowd leaves the domain well before time $T$, so that we do not introduce alien temporal boundary conditions in the HJB equation.

Similarly, a bounded domain $\Omega$ (instead of $\R^2$) must be used in the finite-horizon problem. Boundary conditions in the HJB equation are imposed by setting $\phi=+\infty$ on the boundary cells. 
In such a way reaching the boundary is never optimal for the crowd.

Regarding the FP equation, instead, we refer the reader to \cite{cristiani2017} for a general treatment of the boundary.

\subsection{Notations.}
Let us denote by $(x^1,x^2)$ the two components of the space vector $x$.
The computational domain $\Omega\times[0,T]$ is divided in cells of side $\Dx \times \Dy \times \Dt$, where $\Dx$, $\Dy$ are the space steps and $\Dt$ is the time step. 
Let us assume that $T$ is a multiple of $\Dt$ to avoid rounding, and define $n_T:=\frac{T}{\Dt}\in\mathbb N$. 
Similarly, we choose $\Omega$ as a rectangular domain of size $n^1\Dx\times n^2\Dy$, with $n^1, n^2\in\mathbb N$.
The generic cell is defined as 
$$
C^n_{i,j}:=
\left[x^1_i-\frac{\Dx}{2}, x^1_i+\frac{\Dx}{2}\right) \times 
\left[x^2_j-\frac{\Dy}{2}, x^2_j+\frac{\Dy}{2}\right) 
\times \left[t^n, t^{n+1}\right)
$$ 
where, as usual,
$$
\{x^1_0,\ldots,x^1_{n^1}\},\qquad
\{x^2_0,\ldots,x^2_{n^2}\},\qquad
\{t^0,\ldots,t^{n_T}\}$$
are the grid points along each dimension, respectively.
We also denote by $\r^n_{i,j}$  the approximation of the value $\r\big((x^1_i,x^2_j),t^n\big)$. 
Analogous notation is also used for the other functions $\phi$, $g$, $\ell$, $V$, $\alpha$.
Finally, let us define the set of indices
$$
I^\mathcal S_{i,j}(\alphastar):=\{r,s \ : \ C_{r,s}\in\mathcal S_{i,j}(\alphastar)\}
$$
which will be used to refer to the sensory region \eqref{def:sensoryregion} on the grid.
\subsection{Fokker-Planck equation.}\label{sec:numerics:FP}
Regarding the advection part of the Fokker-Planck equation, we consider the push-forward numerical scheme originally proposed in \cite{piccoli2011ARMA} and then used in \cite{cristiani2017, cristiani2011, cristiani2014book}. 
Although it exhibits a diffusive behavior, it is able to reproduce the main features of pedestrian flow such as merging and splitting.

The scheme is used for both the forward equations, i.e.\ for approximating both $\r$ and $\tilde\r$. Let us focus on $\r$ for notation convenience.

The push-forward scheme is truly two-dimensional and follows the physics of the underlying problem: at each time step, some mass leaves the cell $C_{i,j}$, while other mass coming from neighboring cells enters $C_{i,j}$. A CFL condition of the form
\begin{equation}\label{CFL}
\Dt \ \max_{i,j,n} \{V_{i,j}^n\}\leq \min\{\Dx,\Dy\} 
\end{equation}
is imposed to avoid that any pedestrian mass covers a distance larger than a cell in one time step.

The balance of mass is given by \cite[Sect.\ 5.5.2]{cristiani2014book}
\begin{equation}\label{schema_push_forward}
\r^{n+1}_{i,j} = \frac{1}{\Dx \Dy}
\sum\limits_{r,s\in I^\mathcal S_{i,j}(\alphastar)}
\r^{n}_{i,j}\ \mathcal L\big(C_{i,j}\cap\gamma^n(C_{r,s};\alpha^{*,n}_{i,j},\rho^n)\big), \qquad n=0,\ldots,n_T-1,
\end{equation}
where $\mathcal L(\mathcal K)$ denotes the Lebesgue measure of any subset $\mathcal K\in\R^2$, 
and 
$$
\gamma^n(\zeta;a,\rho^n):=\zeta+V(\zeta,n\Dt;a,\rho^n)\Dt,\qquad \zeta\in\R^2,\quad a\in A
$$
is the discrete-in-time map which pushes the mass forward. 

The scheme \eqref{schema_push_forward} can be written in a more computer-friendly form as follows: 
let us denote by $\delta_{i,j}$ the standard Kronecker delta and by $(\phantom{x})^\pm$ the positive/negative part. Let us also define $X:=\Dt V$ and denote by $X^1, X^2$ are the two components of the vector $X$.
Then, the scheme \eqref{schema_push_forward} can be written as
\begin{equation}\label{schema_push_forward_implementato}
\r^{n+1}_{i,j} =\frac{1}{\Dx\Dy} \sum\limits_{r,s\in I^\mathcal S_{i,j}(\alphastar)} \r^{n}_{r,s} \ \Gamma_{r,s}^{1,n} \ \Gamma_{r,s}^{2,n},
\qquad n=0,\ldots,n_T-1, 
\end{equation}
where
\begin{eqnarray*}
\Gamma_{r,s}^{1,n} := & & 
\left(X_{r,s}^{1,n}\right)^+ \delta_{r,i-1} + 
\left(X_{r,s}^{1,n}\right)^- \delta_{r,i+1} + 
\left(\Dx-|X_{r,s}^{1,n}|\right)\delta_{r,i}, \\ [2mm]
\Gamma_{r,s}^{2,n} := & &
\left(X_{r,s}^{2,n}\right)^+ \delta_{s,j-1} + 
\left(X_{r,s}^{2,n}\right)^- \delta_{s,j+1} + 
\left(\Dy-|X_{r,s}^{2,n}|\right)\delta_{s,j}.
\end{eqnarray*}

Regarding the diffusion part, the Laplacian operator is discretized by a standard 5 points finite difference approximation:
$$
\triangle\rho\big((x^1_i,x^2_j),t^n\big)\approx
\frac{1}{\Dx\Dy} \big(\rho^{n}_{i+1,j}+\rho^{n}_{i-1,j}+\rho^{n}_{i,j+1}+\rho^{n}_{i,j-1}-4\rho^{n}_{i,j}\big).
$$

All together, the final scheme reads as
\begin{equation}\label{schema_push_forward_implementato_CI}
\left\{
\begin{array}{l}
\displaystyle
\r^{n+1}_{i,j} =\frac{1}{\Dx\Dy} \sum\limits_{r,s\in I^\mathcal S_{i,j}(\alphastar)} \r^{n}_{r,s} \ \Gamma_{r,s}^{1,n} \ \Gamma_{r,s}^{2,n}  +\\ [7mm] \qquad\qquad\qquad
+\frac{\sigma\Dt}{\Dx\Dy} \big(\rho^{n}_{i+1,j}+\rho^{n}_{i-1,j}+\rho^{n}_{i,j+1}+\rho^{n}_{i,j-1}-4\rho^{n}_{i,j}\big),
\qquad n=0,\ldots,n_T-1 \\ [5mm]
\displaystyle
\r^0_{i,j} = \frac{1}{\Dx\Dy}\iint\limits_{C_{i,j}}\r_0(y)dy 
\end{array}
\right.
\end{equation}

\subsection{Hamilton-Jacobi-Bellman equation.}\label{sec:numerics:HJB}
For the HJB equation we employ a first-order semi-Lagrangian scheme \cite{cacace2014, falconebook}. For simplicity, in the following we drop the dependency of all variables on $\rtheta$.
The scheme is easily built using backward-in-time first-order finite difference for time derivative
$$
\partial_t\phi\big( (x^1_i, x^2_j),t^n\big) \approx 
\frac{\phi^{n}_{i,j}-\phi^{n-1}_{i,j}}{\Dt},
$$ 
and a discrete partial derivative in space
$$
V\big( (x^1_i, x^2_j),t^n;a\big)\cdot\nabla\phi\big( (x^1_i, x^2_j),t^n\big) \approx
\frac{\phi^n\big(\omega^n_{i,j}(a)\big)-\phi^n_{i,j}}{h},
$$
where $h>0$ is small and
$$
\omega^n_{i,j}(a):=(x^1_i, x^2_j)+h V^n_{i,j}(a)
$$
($\omega^n_{i,j}$ corresponds to the foot of characteristics obtained integrating the controlled advection dynamics by means of a first-order explicit Euler approximation for a time step $h$.)
So we get 
$$
-\frac{\phi^{n}_{i,j}-\phi^{n-1}_{i,j}}{\Dt}=
\min_{a\in A}
\left\{\frac{\phi^n\big(\omega^n_{i,j}(a)\big)-\phi^n_{i,j}}{h}+\ell^n_{i,j}(a)\right\}
+\sigma[\triangle\phi]^n_{i,j}
$$
and then, choosing $h=\Dt$, multiplying by $\Dt$, and discretizing the Laplacian operator as before, we get  
\begin{equation}\label{schema_SL}
\left\{
\begin{array}{l}
\displaystyle
\phi^{n-1}_{i,j} = 
\min\limits_{a\in A}\left\{\phi^n\big(\omega^n_{i,j}(a)\big)+\Dt \ \ell^n_{i,j}(a)\right\}+ \\ [3mm] 
\qquad\quad +\frac{\sigma\Dt}{\Dx\Dy} \left(\phi^{n}_{i+1,j}+\phi^{n}_{i-1,j}+\phi^{n}_{i,j+1}+\phi^{n}_{i,j-1}-4\phi^{n}_{i,j}\right), 
 \qquad n=n_T,n_T-1,\ldots,2,1 \\ [5mm]
\phi^{n_T}_{i,j} = g(x^1_i,x^2_j)
\end{array}
\right. 
\end{equation}
One further step is needed to make the scheme fully discrete: 
in fact, the point $\omega^n_{i,j}$ is not, in general, a grid point. Therefore the value of the function $\phi$ at time $t^n$ at that point must be interpolated using only values of $\phi$ at grid points only. 
We do that using a bilinear interpolation which uses the four vertexes of the cell in which the point falls.
Finally, we approximate the set $A$ with a discrete set $\{a_1,\ldots,a_K$\} constituted by $K$ points, and then we perform the search for the minimum simply comparing all the $K$ possible values.

Once the solution to \eqref{schema_SL} is computed, the optimal control is given by
\begin{equation}\label{synthesis_numerics}
\alpha^{*,n}_{i,j}=\arg\min\limits_{a\in A}\left\{\phi^n\big(\omega^n_{i,j}(a)\big)+\Dt \ \ell^n_{i,j}(a)\right\}.
\end{equation}
\subsection{Mean-field game and convergence to Nash equilibrium.}\label{sec:convergence}
Once the numerical schemes for the FP and the HJB equations are set, they must be duly coupled to solve the forward-backward system \eqref{model:mfg}.
We adopt an iterative scheme which alternates the computation of the solution of the forward equation and the backward equation, starting with the backward one. 
\begin{center}
\begin{tikzcd}
\phi \arrow[to=1-2] 
& \alphastar \arrow[d, "{\text{\tiny FP}}"]\\
\rho^\theta \arrow[from=2-2] \arrow[to=1-1, "{\text{\tiny HJB}}"]
& \tilde\rho
\end{tikzcd}
\end{center}

Another index $k=1,2,\ldots$ is used to count the iterations. 
We denote by $\rho_{(k)}^\theta$ and $\phi_{(k)}$ the solutions to the system after iteration $k$. 
Convergence of the forward-backward system is monitored evaluating the (discrete version of the) distance $E_k$ defined by 
\begin{equation}\label{def:Ek}
E_k:=\|\r_{(k)}^\theta-\r_{(k-1)}^\theta\|_{L^1(\Omega\times[0,T])}.
\end{equation}

It is known that this kind of iterative algorithm can be not convergent to the solution of the problem even if the problem is provably well-posed. 
If this happens, the solution can be stabilized employing the \emph{fictitious play} strategy  \cite{brown1951,cardaliaguet2017}, i.e.\ by passing to the HJB equation a linear combination of the previously computed solutions of the FP equations $\r_{(k)}^\theta,\r_{(k-1)}^\theta,\r_{(k-2)}^\theta,\ldots$. 
In that learning procedure, agents keep memory of the previous choices of the others, and not only of the last ones.

Conversely, if the problem has multiple solutions (i.e.\ multiple Nash equilibria), it is likely that the algorithm starts oscillating between two of them and never stabilizes. In the latest case, the distance $E_k$ converges to a nonzero value or does not converge at all.

\medskip

Summarizing, the complete algorithm to solve equations \eqref{model:main_forward}-\eqref{model:mfg} reads as
\begin{enumerate}
	\item Set $n=0$ and initialize $\rho^0$;
	\item Iterate the solution to FP and HJB in \eqref{model:mfg} with schemes \eqref{schema_push_forward_implementato_CI} and \eqref{schema_SL} until $E_k$ tends to 0 or at least stabilizes \REV{(otherwise a fixed maximal number of iterations is adopted)};
	\item Use $\alpha^{*,n}$ computed in \eqref{synthesis_numerics} to move forward the main FP equation \eqref{model:main_forward} again with scheme \eqref{schema_push_forward_implementato_CI} from time step $n$ to time step $n+1$;
	\item Go to step 2 with $n\leftarrow n+1$, until $n_T$ is reached.
\end{enumerate}

\medskip

\REV{\begin{rem}\label{rem:numericaltrick}
	Actually it is not really needed to run the algorithm until the final time $T$ (i.e.\ for $n_T$ time steps).
	In fact, the main forward equation \eqref{model:main_forward} can be stopped at time $T-\theta$ since, at that time, we already know the complete solution thanks to prediction until time $T$ obtained by solving \eqref{model:mfg}. This allows to save a certain amount of computational time.
\end{rem}}

\begin{rem}
	In the case of the MFG with limited predictive abilities (i.e.\ $0<\theta<T$), several possibilities can appear: 
	it is possible that, for a fixed $\theta$,  the forward-backward system is well-posed for some time instants and ill-posed for other time instants. 
	Alternatively, it is possible that for a fixed time instant $t$,  the model is well-posed for some $\theta$'s and ill-posed for other $\theta$'s. 
	Numerical tests will give some clues on this point.
\end{rem}

\section{Numerical simulations.}\label{sec:tests}	
In this section we present five numerical tests for both finite-horizon and minimum-time problems.
We restrict ourselves to observe the qualitative behaviour of the crowd, leaving more quantitative and practical aspects to a future research.
In all tests we set $A=B(0,1)$ (discretized by $K=32$ points, all on the boundary of the ball), $\Omega=[0,1]^2$ with $n^1=n^2=50$ ($\Dx=\Dy=0.02$), $R_0=0.01$, $R=0.06$.
\subsection{Finite-horizon problem.}
Here we focus on the finite-horizon scenario with and without diffusion term. 
Test 1 investigates the connection between uniqueness of the Nash equilibrium and convergence of the numerical algorithm. 
Test 2 highlights the effect of the key parameter $\theta$ of the proposed model on the pedestrian dynamics.
\subsubsection{Test 1: Well-posedness.}
In this test we choose $T=0.5$ (with $n_T=600$), $\sigma=0.05$, $\ell(\rho)=3 \rho$, and $g(x^1,x^2)=\sqrt{(x^1-0.5)^2+(x^2-0.5)^2}$. 
At initial time, pedestrians are arranged in a square of side 0.1 in the left-bottom corner of the domain, see Fig.\ \ref{fig:rho0,ell,g}(a).
The repulsion parameter is $C_\textsc{rep}=0$, therefore we get $\Vi=0$. This means that the coupling between FP and HJB is all in the running cost $\ell$.

In this case, the structure of the MFG system (obtained with $\theta=T$), falls into the framework in which existence and uniqueness of the solution have been proved, see \cite{lasry2007}.

As expected, the crowd starts moving towards the center of the domain, diffusing around its barycenter, due to the diffusion term and the running cost itself.
In this test we do not show the crowd behaviour, which is rather trivial, instead we focus on the convergence of the algorithm. 
Hereafter, we will say that the algorithm converges if  $E_k\to 0$ as $k\to +\infty$, see \eqref{def:Ek}.

We have performed several simulations for different values of $\theta$, ranging from $\theta=0$ to $\theta=T$. 
Numerical evidence suggests the existence of a threshold value $\overline{\theta}$ for the parameter $\theta$ in order to guarantee convergence. 
In particular, we observe that our algorithm converges, \emph{for every time step $n$}, only for small values of $\theta$ (in our test $ \theta <\overline{\theta}= 0.23$). 
Otherwise, the algorithm does not converge, i.e.\ $E_k\to \bar E>0$ as $k\to +\infty$.

Employing the \REV{stabilizing} fictitious play strategy (see Sect.\ \ref{sec:convergence}), instead, we recover convergence at each time step and for any $\theta$. 
This is in agreement with the theoretical results which assure well-posedness of the theoretical problem as well as the convergence of the numerical algorithm, at least for the MFG case ($\theta=T$).

\REV{In conclusion, \emph{numerics seem to suggest that the problem is well-posed for any $\theta$}, but the algorithm becomes less and less stable as $\theta$ increases, therefore some stabilization strategy is needed to get numerical convergence.}


\REV{\paragraph{Computational time.} Since in this test the algorithm converges in a finite number of steps, it is possible to compute the CPU time with no ambiguity. We run the numerical code on a laptop equipped with an Intel Core i7-1060NG7 processor and 16 GB RAM. The code is serial. 
\begin{itemize}[leftmargin=*]
\item For $\theta=0$ we do not need forward-backward iterations and the computation until final time $T$ takes 7 m 50 s.
\item For $\theta=0.08$ we need 7 iterations to complete every forward-backward step \eqref{model:mfg}. We stop the main forward equation \eqref{model:main_forward} at time $t=T-\theta=0.42$, see Remark \ref{rem:numericaltrick}. Computation takes 5 h.
\item For $\theta=0.2$ we need 17 iterations to complete every forward-backward step \eqref{model:mfg}. We stop the main forward equation \eqref{model:main_forward} at time $t=T-\theta=0.3$, see Remark \ref{rem:numericaltrick}. Computation takes 9 h 55 m.
\item For $\theta=T$ we need 49 iterations to complete the only forward-backward step \eqref{model:mfg}. The main forward equation \eqref{model:main_forward} is not solved at all. Computation takes 4 m 54 s.
\end{itemize}
We can see that the feature of limited predictive abilities enormously increases the computational effort, because of the large number of MFGs \eqref{model:mfg} to be solved.
}
\subsubsection{Test 2: The role of parameter $\theta$.}
In this test we show the behaviour of the crowd for different values of $\theta$. 
We choose  $T=1$ (with $n_T=200$), $\sigma=0$, $\ell (x^1)=-2x^1+3$, and $g(x^1,x^2)=\sqrt{(x^1-0.5)^2+(x^2-0.5)^2}$, while the initial distribution of pedestrians is as in Test 1, see Fig.\ \ref{fig:rho0,ell,g}.
\begin{figure}[h!]
	\centering
	\subfigure[][Initial configuration $\r_0$]{\includegraphics[width=3.5cm, height=3.3cm]{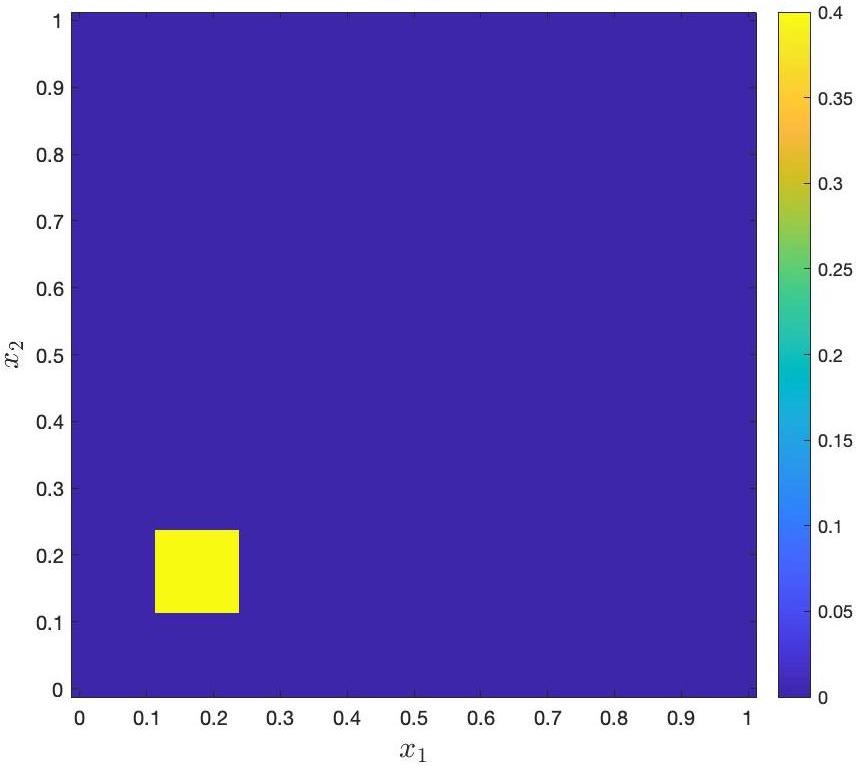}}
	\hspace{0 mm}
	\subfigure[][Running cost $\ell$]{\includegraphics[width=3.5cm, height=3.3cm]{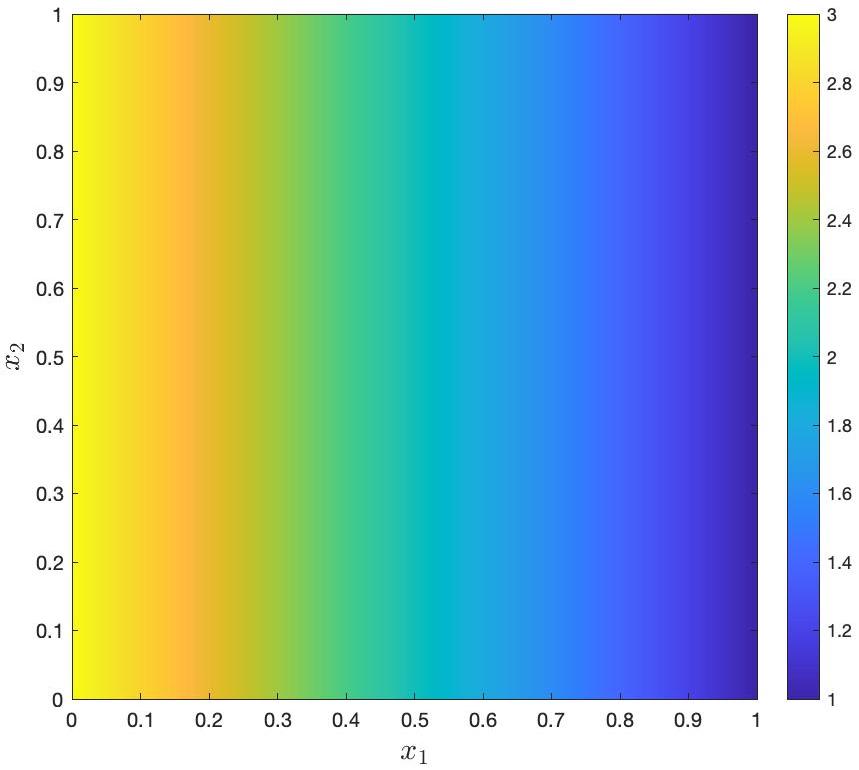}}
	\hspace{0 mm}
	\subfigure[][Final cost $g$]{\includegraphics[width=3.5cm, height=3.3cm]{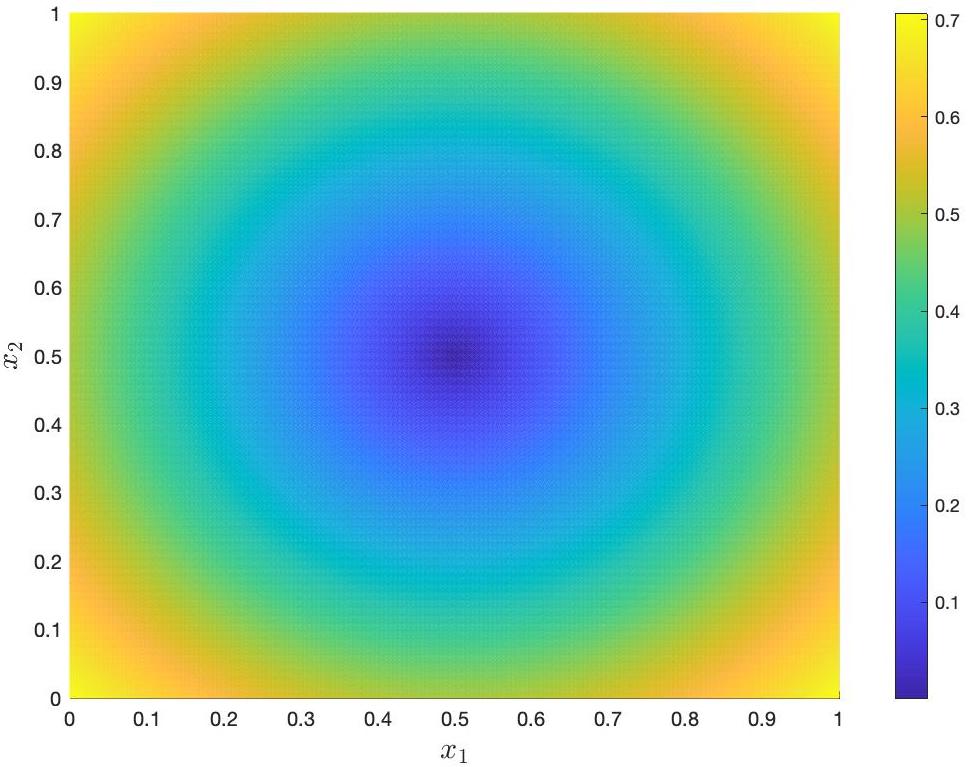}}
	\caption{Test 2: Given functions}
	\label{fig:rho0,ell,g}
\end{figure}
The repulsion parameter is $C_\textsc{rep}=6$. 
Note that here $\ell$ does not depend on $\rho$. This means that, conversely to Test 1, the coupling between FP and HJB is all in the interaction velocity $\Vi$.
 
\begin{figure}[h!]
\centering
\raisebox{40pt}{\parbox[b]{.12\textwidth}{$t=0.12$}}
\subfigure[][$\theta=0$]{\includegraphics[width=3.0cm, height=3.0cm]{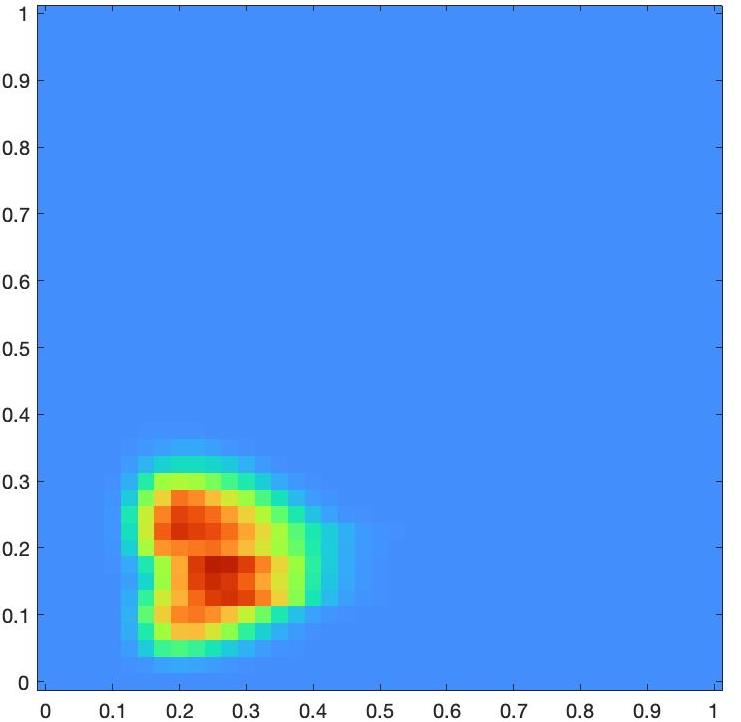}}
\hspace{0 mm}
\subfigure[][$\theta=0.25$]{\includegraphics[width=3.0cm, height=3.0cm]{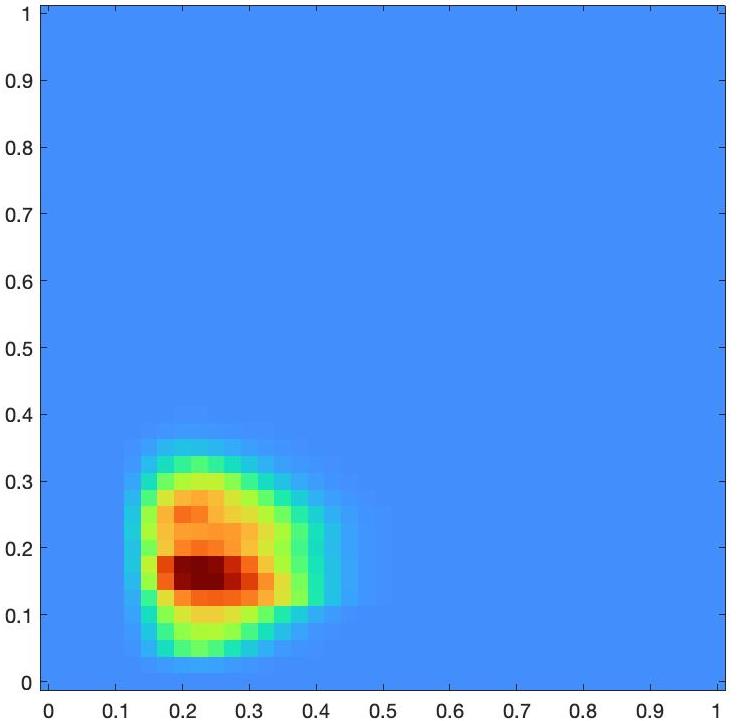}}
\hspace{0 mm}
\subfigure[][$\theta=1$]{\includegraphics[width=3.3cm, height=3.0cm]{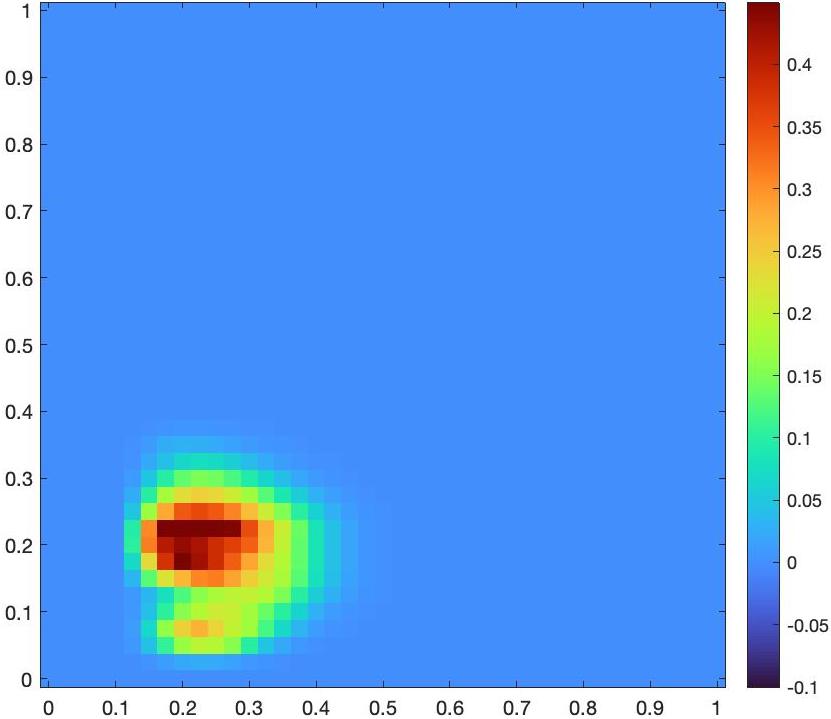}}\\
\vspace{0 mm}
\raisebox{40pt}{\parbox[b]{.12\textwidth}{$t=0.55$}}
\subfigure[][$\theta=0$]{\includegraphics[width=3.0cm, height=3.0cm]{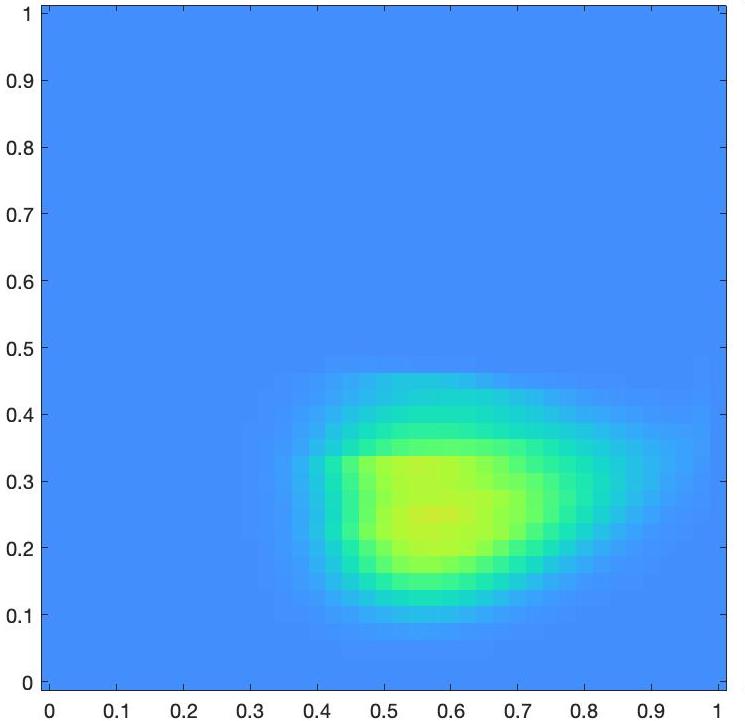}}
\hspace{0 mm}
\subfigure[][$\theta=0.25$]{\includegraphics[width=3.0cm, height=3.0cm]{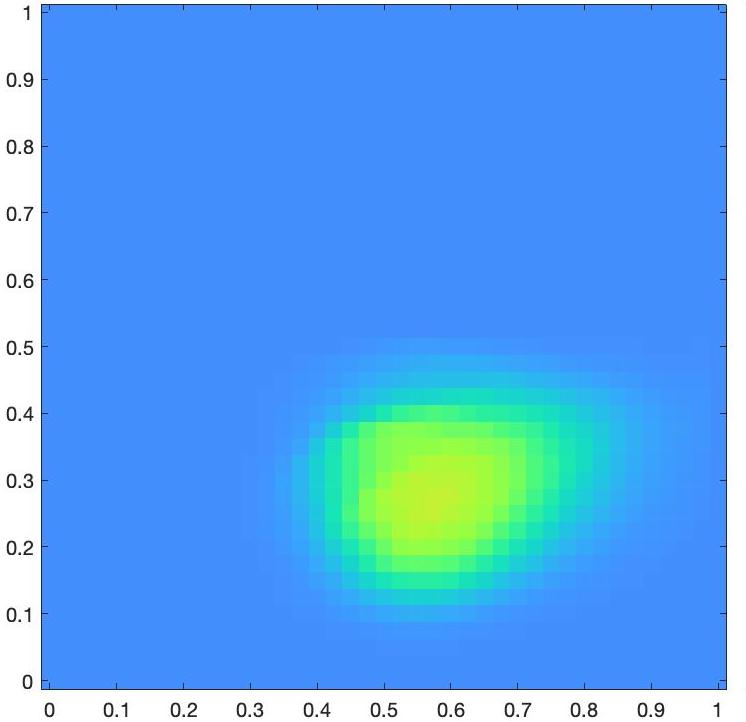}}
\hspace{0 mm}
\subfigure[][$\theta=1$]{\includegraphics[width=3.3cm, height=3.0cm]{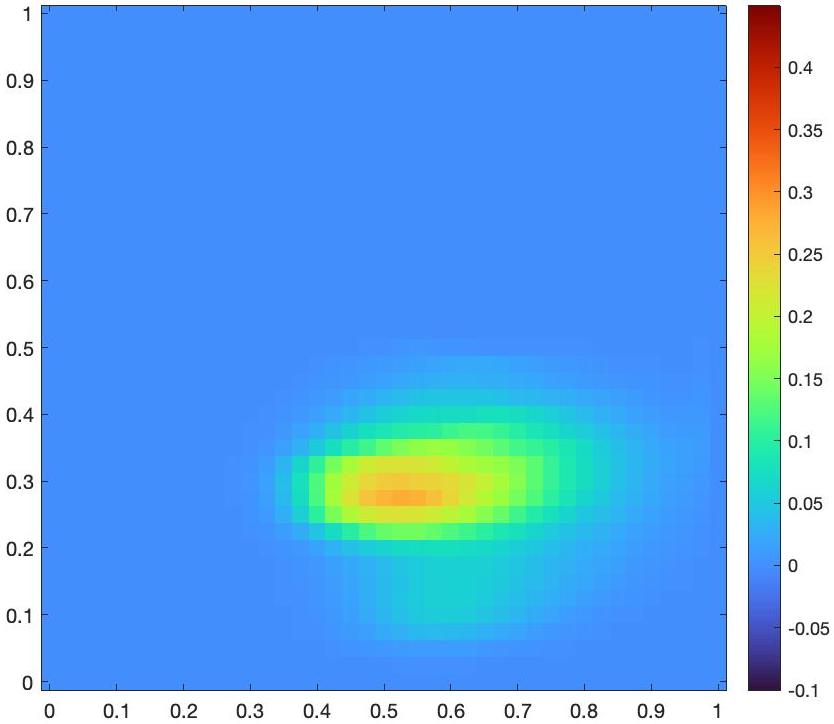}}\\
\vspace{0 mm}
\raisebox{40pt}{\parbox[b]{.12\textwidth}{$t=0.75$}}
\subfigure[][$\theta=0$]{\includegraphics[width=3.0cm, height=3.0cm]{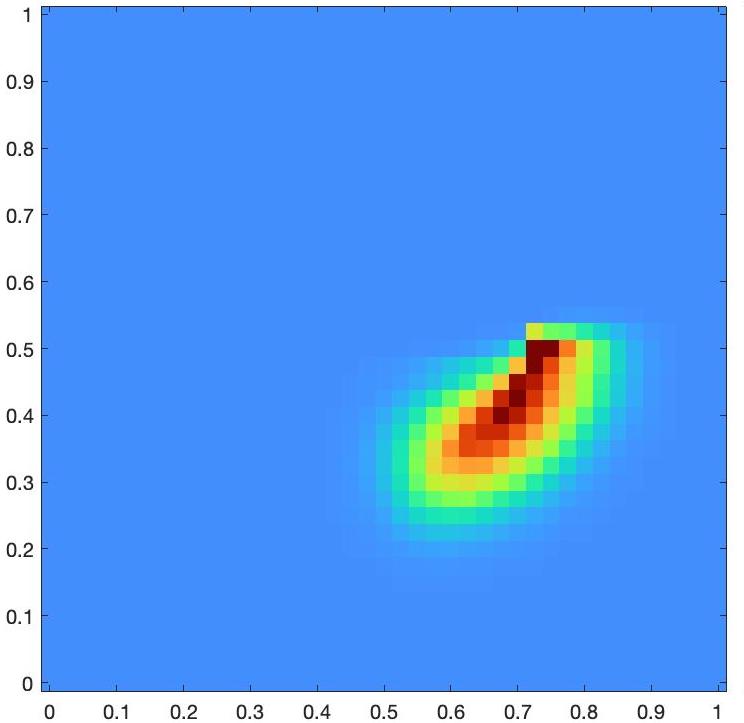}}
\hspace{0 mm}
\subfigure[][$\theta=0.25$]{\includegraphics[width=3.0cm, height=3.0cm]{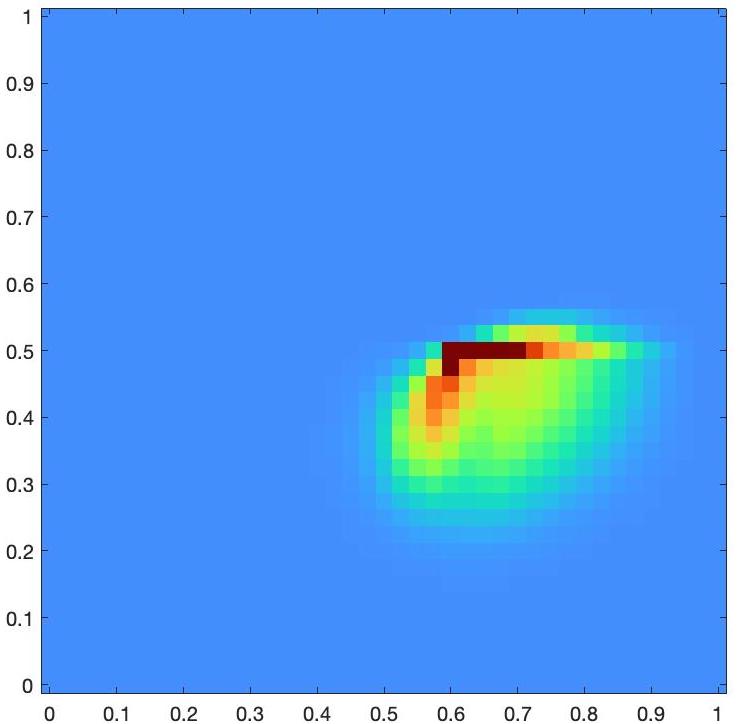}}
\hspace{0 mm}
\subfigure[][$\theta=1$]{\includegraphics[width=3.3cm, height=3.0cm]{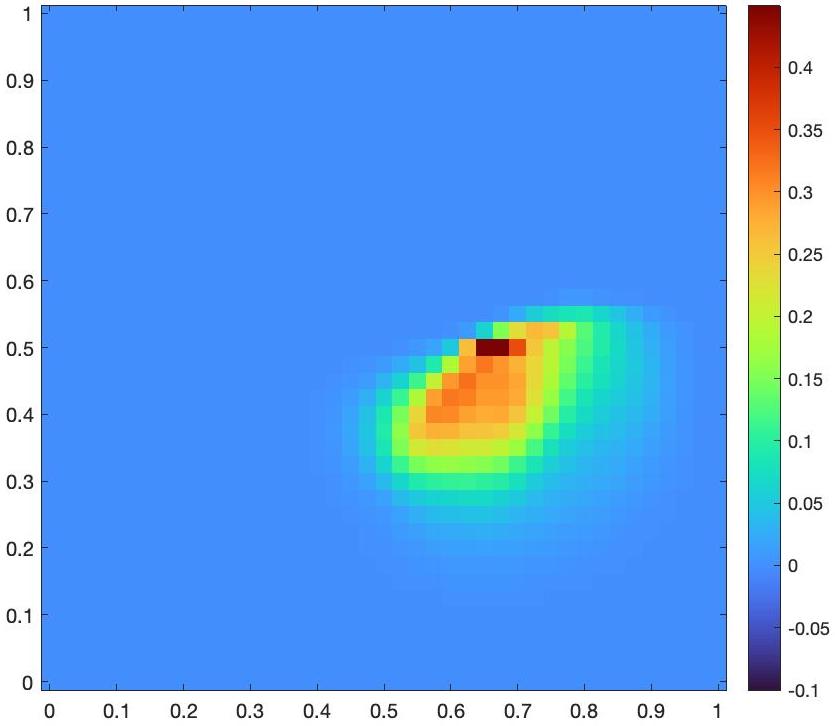}}
\caption{Test 2: Density evolution for three different values of $\theta$}
\label{fig:test2}
\end{figure}

Figure \ref{fig:test2} shows different frames of three representative simulations performed varying the value of $\theta$. 
We compare the case of 
$\theta=0$ (Hughes's type model with no prediction ability), 
$\theta=0.25$, and 
$\theta=T$ (standard MFG with complete prediction ability), 
at three time steps of the simulations. 
%
Numerical results reproduce the expected behavior. In fact, regardless of the value of $\theta$, the crowd moves from the left to the right, in order to reach regions where the running cost has lower values, but, when the final time approaches, the group turns to reach the center of the domain where $g$ is minimal. 

The effect of the time window for the prediction reflects in different strategies adopted by pedestrians to reach the center at the final time:
as $\theta$ increases pedestrians anticipate the change of direction towards the center. In particular, from the very beginning a larger number of pedestrians move from the bottom of the domain toward the center, since this will give them an advantage later on, when they will approach the center of the domain.
Indeed, the possibility to forecast the position of the others allows people to understand that approaching too much the right side (to save the running cost) will create a difficulty to reach the center later on (to save the final cost), because of the large mass of people one will find in front along the way back. 


\subsection{Minimum-time problem.}
Here we consider two different settings in the context of minimum-time problems. 
In Test 3 \& 4 we consider a square room with two exits, and an initial mass of people arranged in two separated groups. As we did in Test 1 \& 2, we first investigate convergence properties of the proposed  algorithm and then the effect of $\theta$ on the dynamics.
Finally, in Test 5, we consider the case of a single time-dependent exit.

\subsubsection{Test 3: Well-posedness.}
\begin{figure}[h!]
\centering
\includegraphics[width=4cm, height=3.6cm]{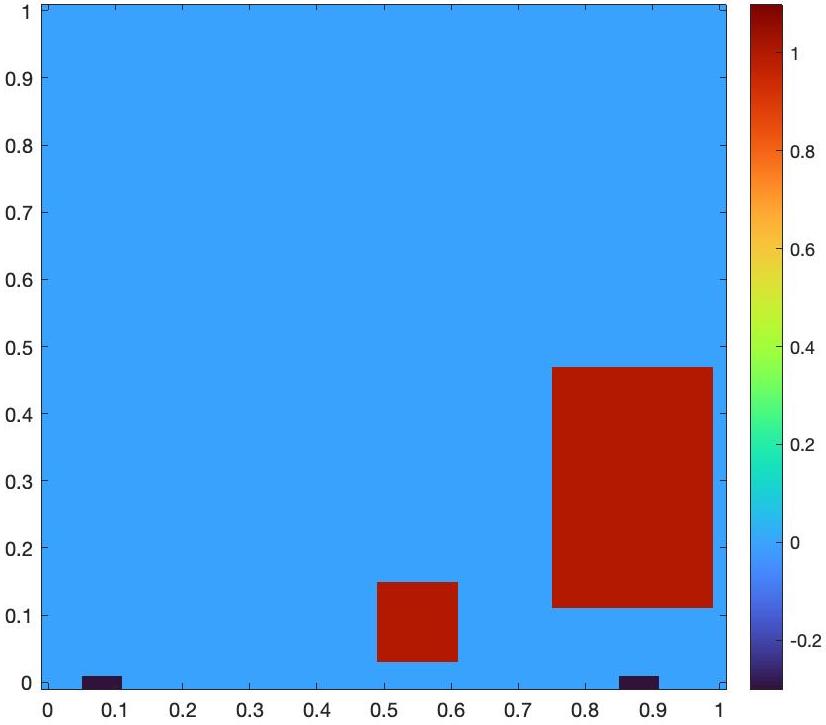}
\caption{Tests 3 \& 4: Initial pedestrian configuration $\r_0$ and exit location along the bottom side}
\label{fig:test4_initial}
\end{figure}

In this test we choose $T=1.5$ ($n_T=200$), $\sigma=0$, $C_\textsc{rep}=8$, and $R=0.06$. 
The initial configuration is constituted by two separated groups located at the bottom and near the right side of the room, respectively, see Fig.\ \ref{fig:test4_initial}. 
Two exits are placed on the bottom side of the domain, one on the left and one the right, respectively.

As already observed in Remark \ref{rem:theroeticalinsights}, uniqueness of the solution for our model is not guaranteed.
As in Test 1,  we observed the behaviour of the distance $E_k$ for every time step $n$ and for every $\theta\in[0,T]$.
The results we have obtained lead to some possible conclusions: first of all, regardless of the value of $\theta>0$, the distance $E_k$ converges to zero only for some time steps of the simulation (sometimes none of them), and this is true even if the fictitious play strategy is employed. 
This suggests that the solution is not unique. 
Moreover, since a different forward-backward system is solved at every time, at every time a different batch of multiple solutions could appear; therefore, a tree-like model of possible outcomes arises. 
This is not true for the standard MFG, where there is only one forward-backward system to be solved, therefore multiple solutions, if any, appear only once and refer to the whole dynamics.  

Secondly, the number of time instants in which convergence occurs decreases as $\theta$ increases. This confirms what we have already observed in Test 1, namely that there is a sort of ``degree of ill-posedness'' which increases together with $\theta$.

In conclusion, \emph{numerics seem to suggest that the problem is ill-posed for any $\theta>0$} and the algorithm becomes more and more unstable as $\theta$ increases. 
Moreover, for intermediate $\theta$'s (i.e.\ $0<\theta<T$), one can combine along the temporal line all the possible solutions obtained at any fixed time.\footnote{This curiously resembles the theory of multiverse in physics.} 
\subsubsection{Test 4: The role of parameter $\theta$.}
\begin{figure}[h!]
\centering
\raisebox{40pt}{\parbox[b]{.12\textwidth}{$t=0.06$}}
\subfigure[][$\theta=0$]{\includegraphics[width=3.0cm, height=3.0cm]{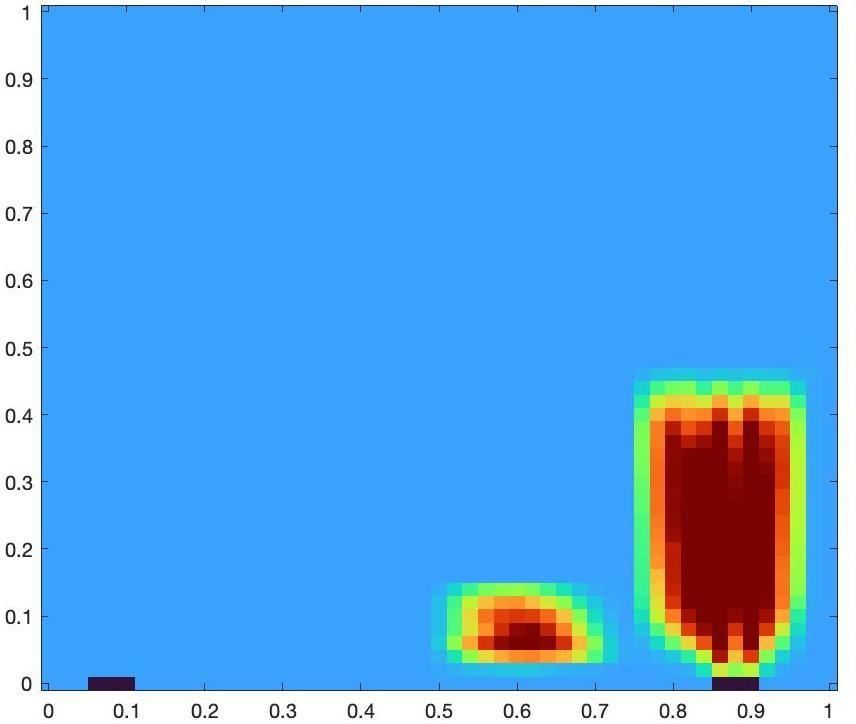}}
\hspace{0 mm}
\subfigure[][$\theta=0.15$]{\includegraphics[width=3.0cm, height=3.0cm]{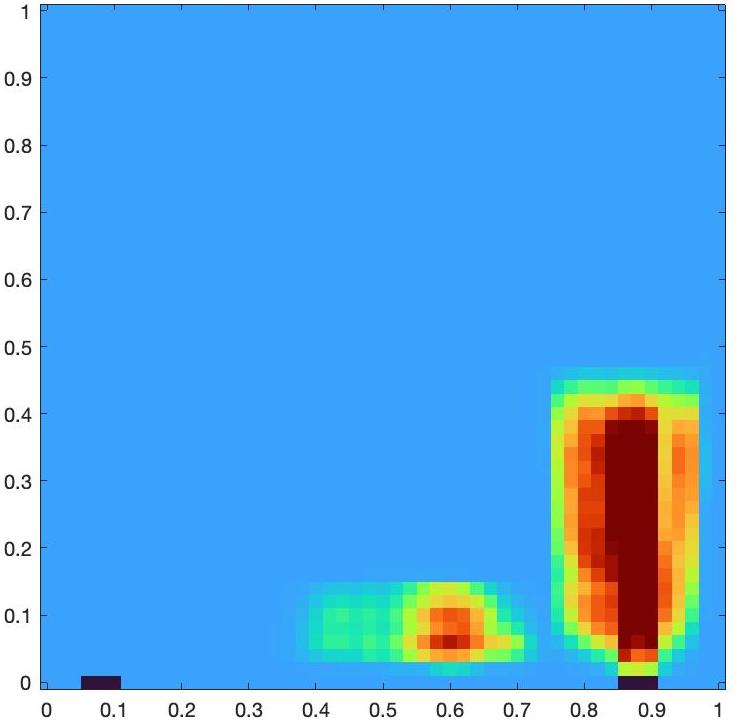}}
\hspace{0 mm}
\subfigure[][$\theta=0.75$]{\includegraphics[width=3.3cm, height=3.0cm]{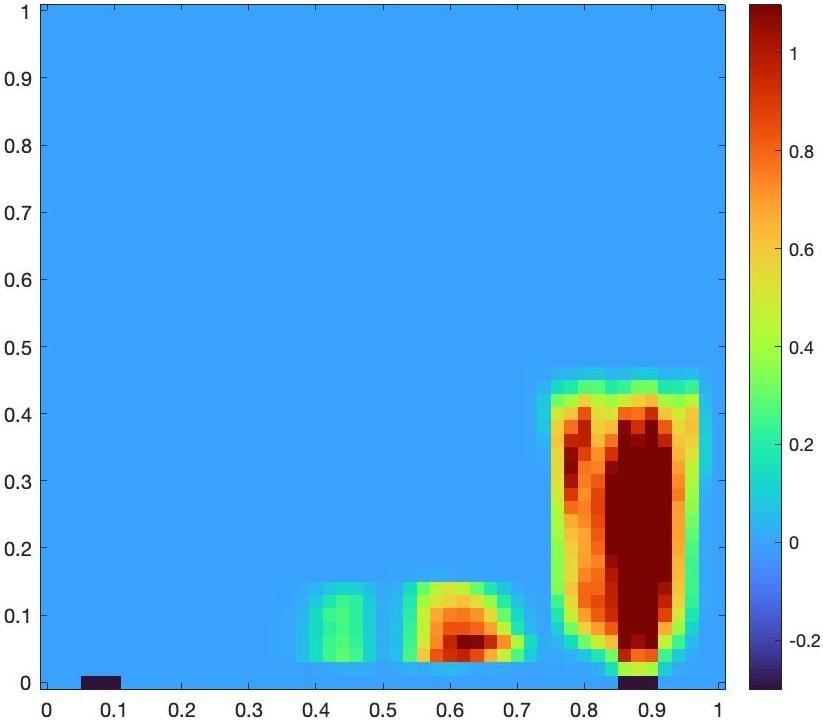}}\\
\vspace{0 mm}
\raisebox{40pt}{\parbox[b]{.12\textwidth}{$t=0.22$}}
\subfigure[][$\theta=0$]{\includegraphics[width=3.0cm, height=3.0cm]{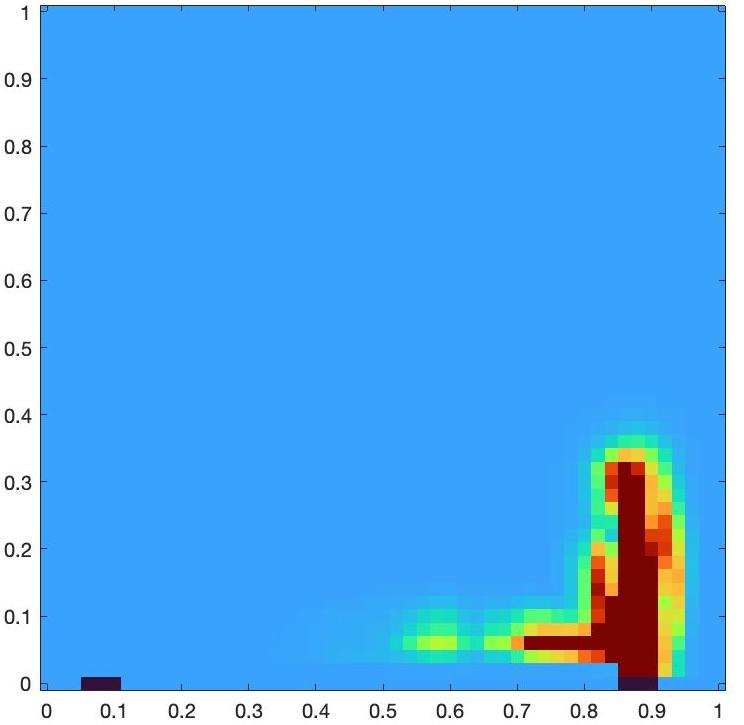}}
\hspace{0 mm}
\subfigure[][$\theta=0.15$]{\includegraphics[width=3.0cm, height=3.0cm]{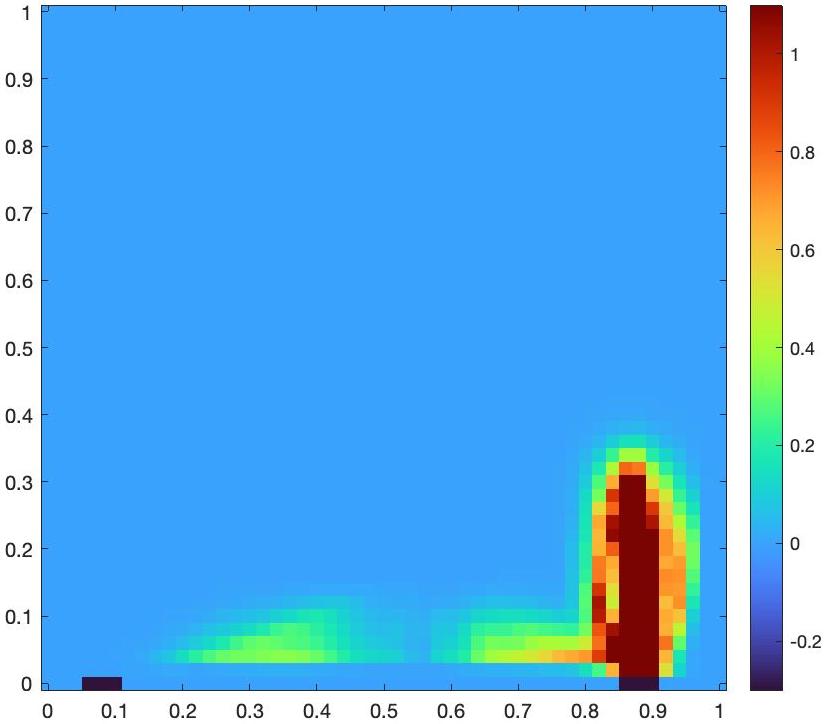}}
\hspace{0 mm}
\subfigure[][$\theta=0.75$]{\includegraphics[width=3.3cm, height=3.0cm]{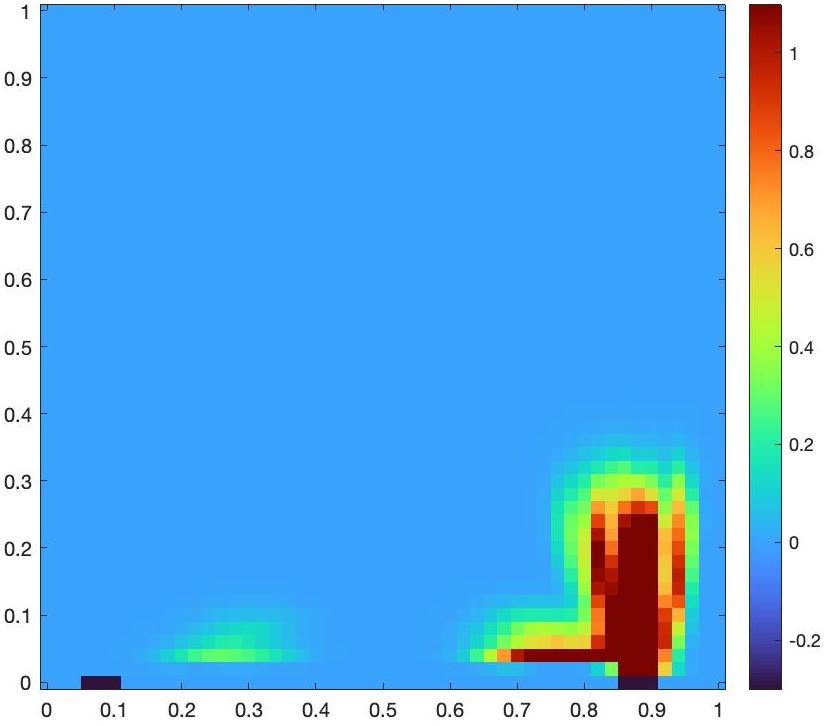}}\\
\vspace{0 mm}
\raisebox{40pt}{\parbox[b]{.12\textwidth}{$t=0.43$}}
\subfigure[][$\theta=0$]{\includegraphics[width=3.0cm, height=3.0cm]{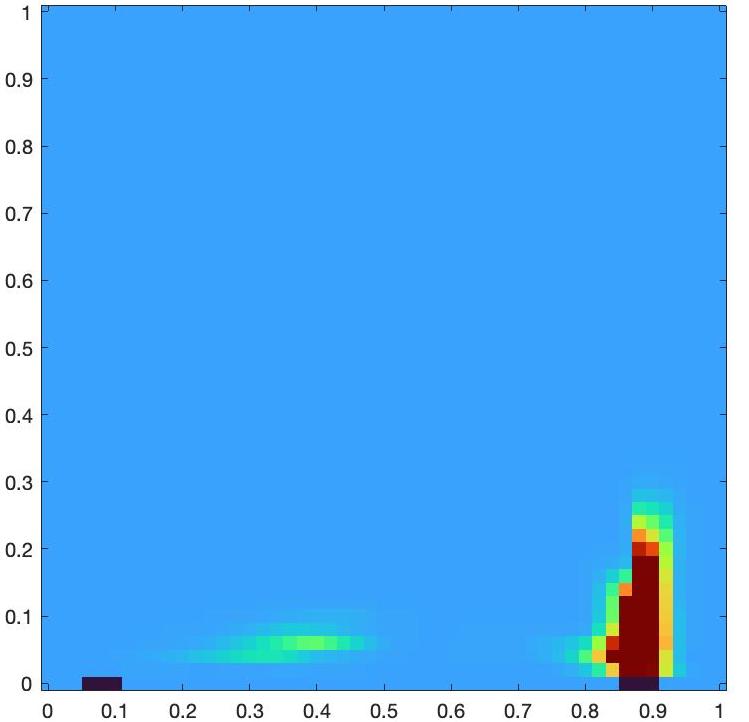}}
\hspace{0 mm}
\subfigure[][$\theta=0.15$]{\includegraphics[width=3.0cm, height=3.0cm]{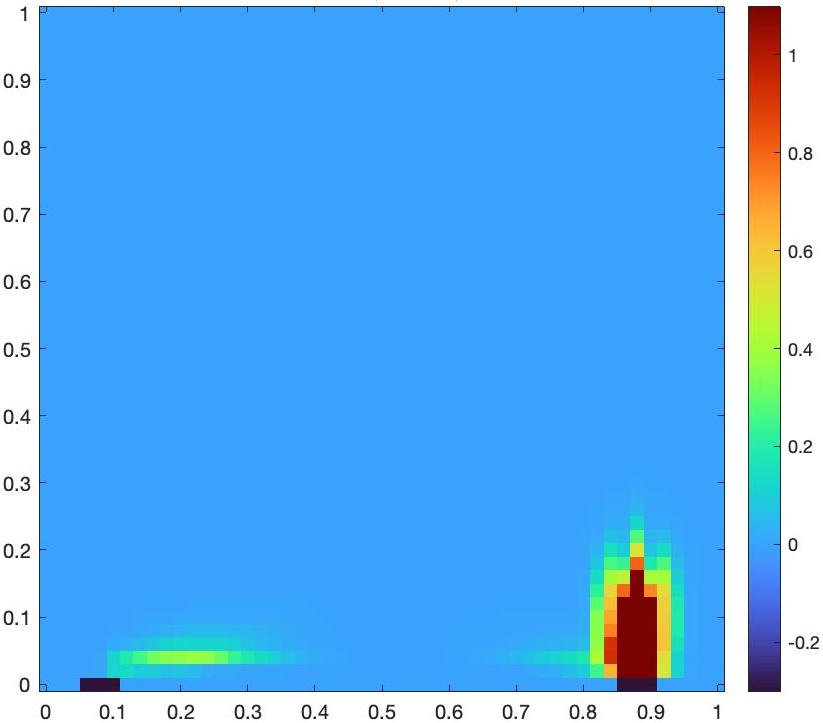}}
\hspace{0 mm}
\subfigure[][$\theta=0.75$]{\includegraphics[width=3.3cm, height=3.0cm]{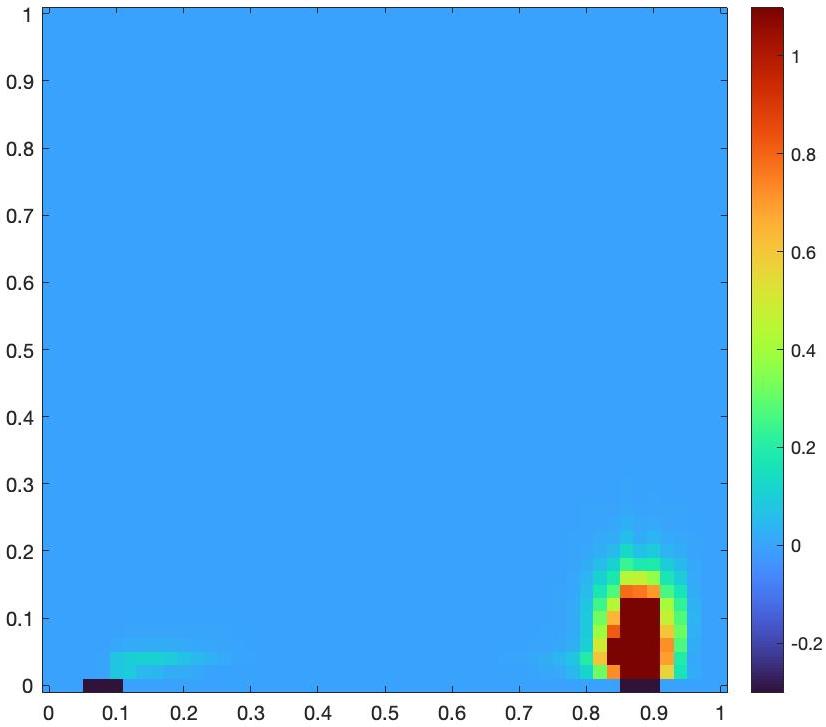}}
\caption{Test 4: Density evolution for three different values of $\theta$}
\label{fig:test4}
\end{figure}

Keeping the same scenario and parameters of Test 3, we now focus on the effect of $\theta$ on the density evolution.
To this end, we run three simulations with $\theta=0, 0.15, 0.75$, see Fig.\ \ref{fig:test4}.
In the case $\theta=0$, pedestrians react to current crowd distribution only and then they are not able to predict the formation of a huge congestion near the rightmost exit. (a) As a consequence, all people move toward that exit, which is the closest one. (d,g) After a while, however, some people understand that the rightmost exit is not convenient and come back, heading toward the leftmost exit, eventually reaching it.

Here the introduction of predictive abilities has a clear effect: some people understand \emph{in advance} that heading left is the best option. 
(b) In the case $\theta=0.15$, some people belonging to the smaller group move leftward from the very beginning, (e,h) while others initially move rightward and U-turn only after a while. Separation occurs gradually as time goes on and people predict the extent of the congestion.
(c,f,i) In the case $\theta=0.75$, instead, people are able to forecast the whole dynamics since the initial time. Therefore the split occurs only at the beginning and no one changes its mind afterwards.
 
(g,h,i) The different behaviour is also visible in the last frame: the larger $\theta$, the shorter the exit time of people moving to the left. This certifies the higher degree of optimality of the dynamics with prediction.

%
%
%
\subsubsection{Test 5: Time-dependent exit.}
\begin{figure}[h!]
	\centering
     \subfigure[]{\includegraphics[width=4cm, height=3.6cm]{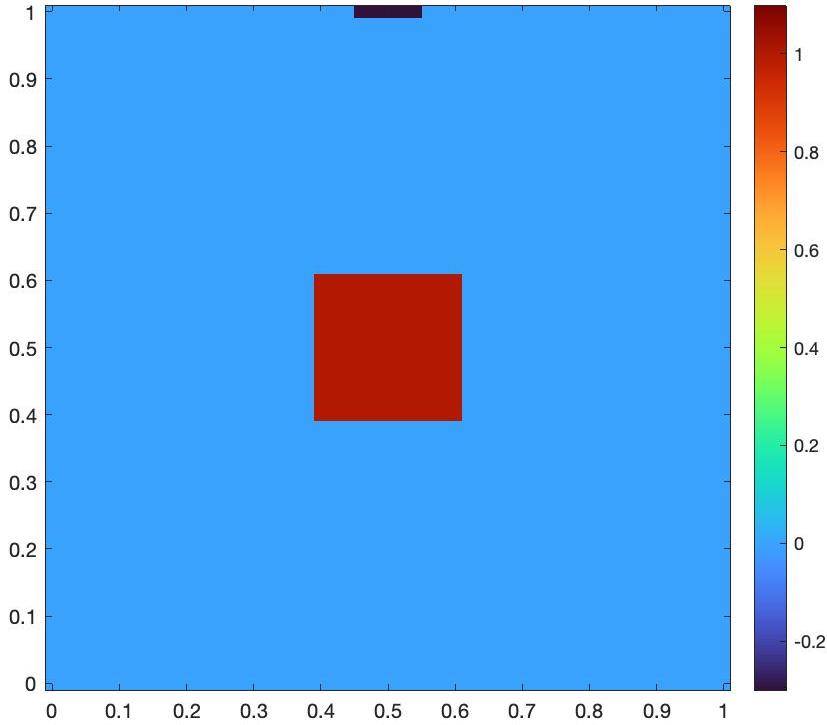}}
	\hspace{10 mm}
	\subfigure[]{
		\begin{overpic}[width=4cm, height=3.6cm]{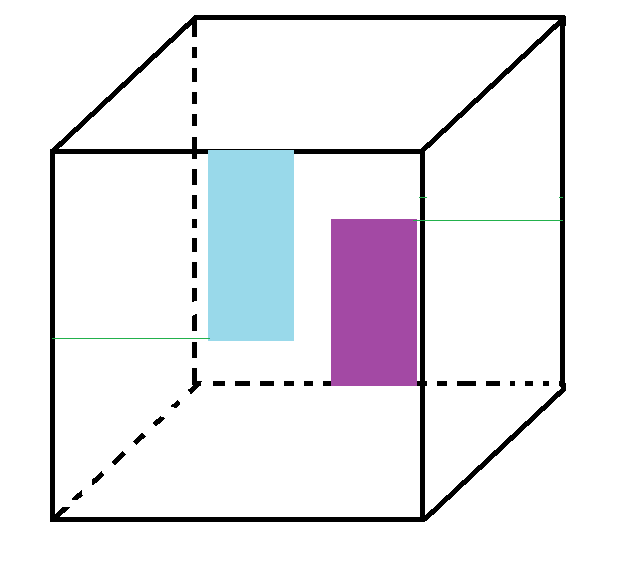}
			\put(76,32){\tiny top} \put(35,10){\tiny bottom}
			\put(55,-3){$x^1$} \put(15,20){$x^2$} \put(0,60){$t$} \put(0,35){$\bar t$}\put(95,52){$\bar t$}
			\put(57,40){$\mathcal T$} \put(37,50){$\mathcal T$}
		\end{overpic}
	}
	\caption{Test 5: (a) Initial configuration $\r_0$. Only the exit along the top side is now open. (b) Target $\mathcal T$ in the space-time: at time $\bar t$ it switches from top side to bottom side}
	\label{fig:test5_initial}
\end{figure}
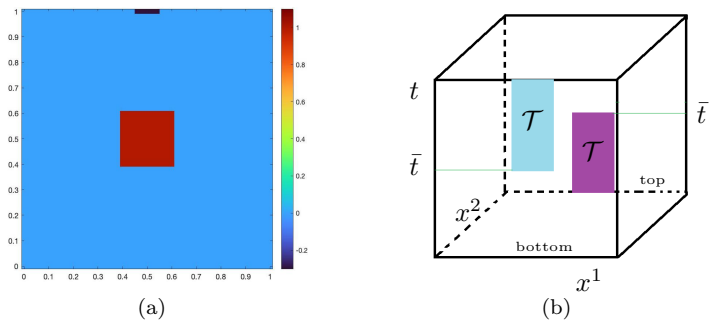

In this test we consider a time-dependent target. 
We choose $T=2.5$ ($n_T=200$), $\sigma=0$, $C_\textsc{rep}=8$. 
The crowd is initially arranged in the center of the domain and the exit is located on the top side of the domain, see Fig.\ \ref{fig:test5_initial}(a). 

At a certain time $\bar{t}$, the exit closes and another exit, in the lower side, opens.
We assume that pedestrians are able to forecast the change of the exit only $\Theta$ time units in advance, where $\Theta>0$ is a novel parameter which combines in a nontrivial manner with $\theta$.

The dependence on time of the exit -- i.e.\ the position of the target $\mathcal T$ -- is obtained easily, since the problem is already set in space-time: $\mathcal T\times[0,T]$ is just a nonconnected set which is located on the top of the domain for $t\in[0,\bar t]$ and on the bottom of the domain for $t\in(\bar t,T]$, see Fig.\ \ref{fig:test5_initial}(b).

The ability of forecasting the change of the exit location is instead more tricky. 
This is obtained by passing to the HJB equation a potentially different target $\mathcal T=\mathcal T(s)$ at every time $s$, i.e.\ every time a forward-backward must be solved. 
Until time $s=\bar t-\Theta$, the target $\mathcal T(s)$ describes an exit fixed at the top of the domain for all times $t\in[0,T]$. 
In this way pedestrians will solve their minimum-time problem believing that the exit will stay on the top forever. 
After time $s=\bar t-\Theta$, the target $\mathcal T(s)$ changes and describes the exit as it is really found, i.e.\ at the top for $t\in[0,\bar t]$ and on the bottom elsewhen. 

\begin{figure}[h!]
\centering
\raisebox{40pt}{\parbox[b]{.12\textwidth}{$t=0.12$}}
\subfigure[][$\theta=0$]{\includegraphics[width=3.0cm, height=3.0cm]{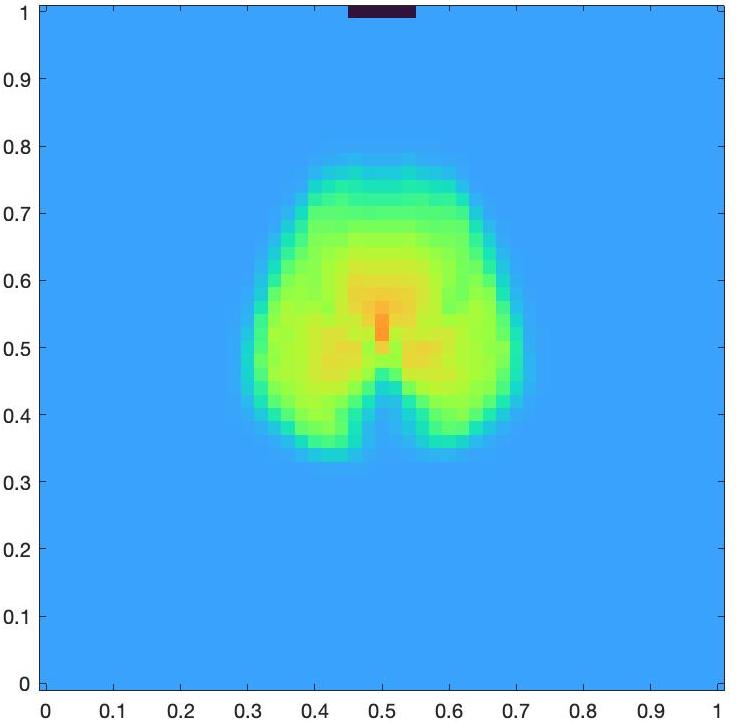}}
\hspace{0 mm}
\subfigure[][$\theta=0.25$]{\includegraphics[width=3.0cm, height=3.0cm]{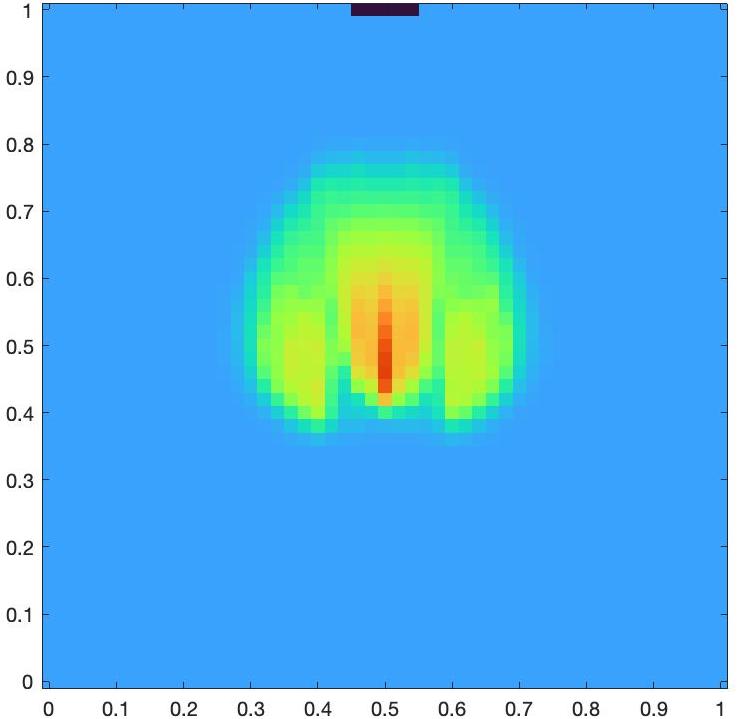}}
\hspace{0 mm}
\subfigure[][$\theta=2.5$]{\includegraphics[width=3.3cm, height=3.0cm]{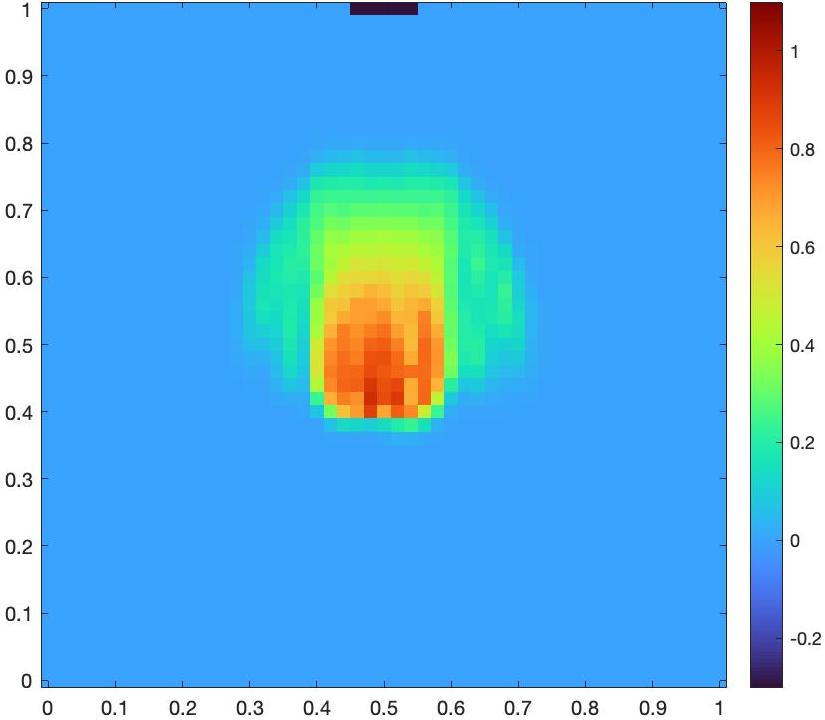}}\\
\vspace{0 mm}
\raisebox{40pt}{\parbox[b]{.12\textwidth}{$t=0.26$}}
\subfigure[][$\theta=0$]{\includegraphics[width=3.0cm, height=3.0cm]{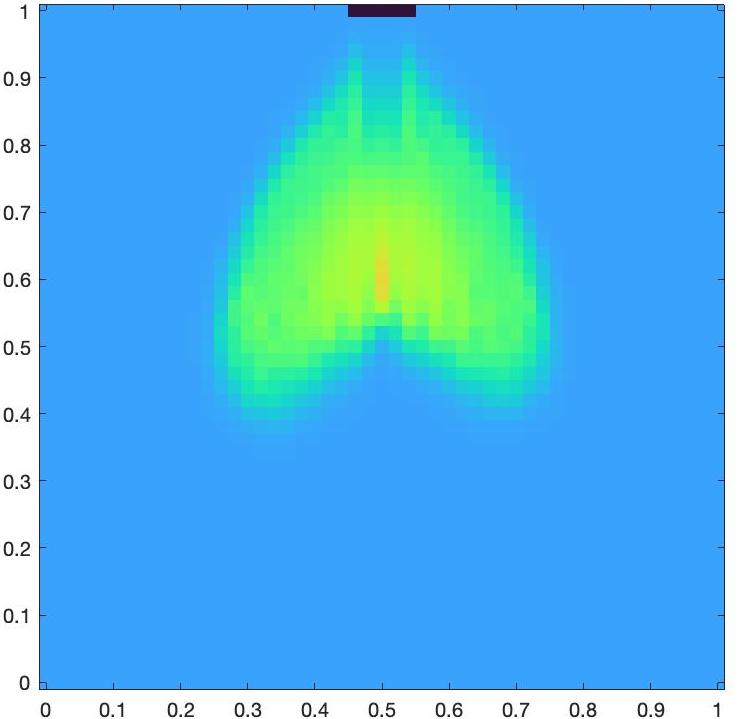}}
\hspace{0 mm}
\subfigure[][$\theta=0.25$]{\includegraphics[width=3.0cm, height=3.0cm]{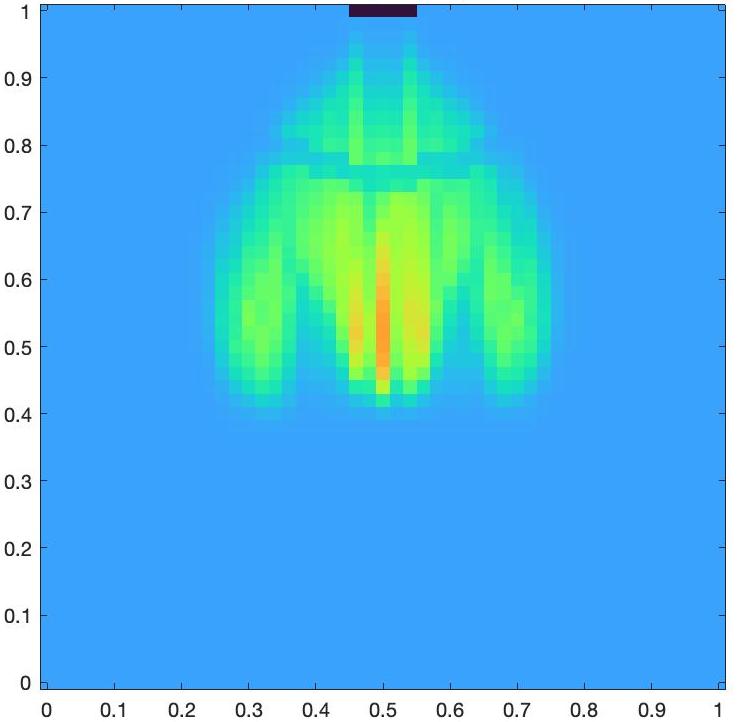}}
\hspace{0 mm}
\subfigure[][$\theta=2.5$]{\includegraphics[width=3.3cm, height=3.0cm]{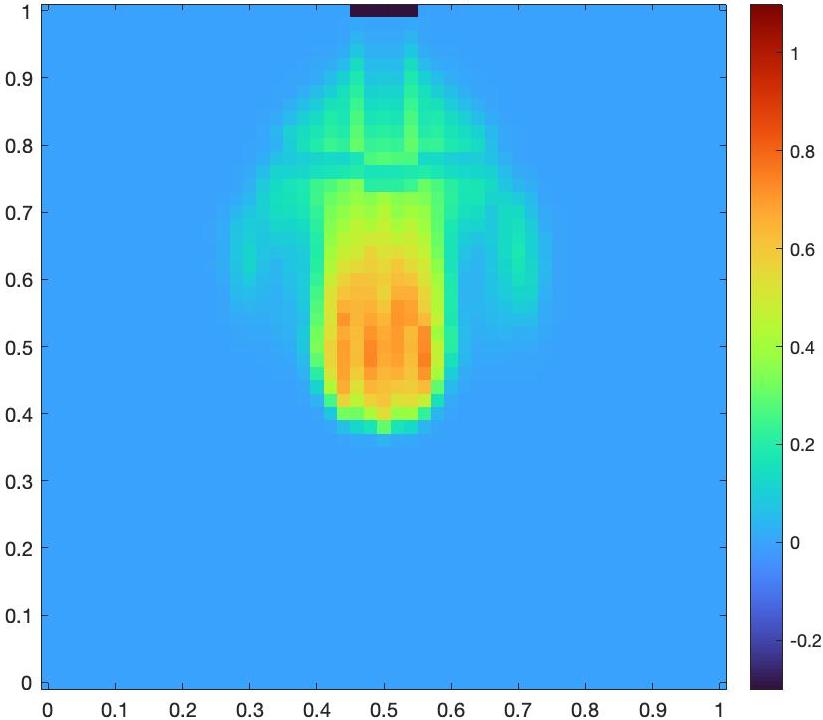}}\\
\vspace{0 mm}
\raisebox{40pt}{\parbox[b]{.12\textwidth}{$t=0.67$}}
\subfigure[][$\theta=0$]{\includegraphics[width=3.0cm, height=3.0cm]{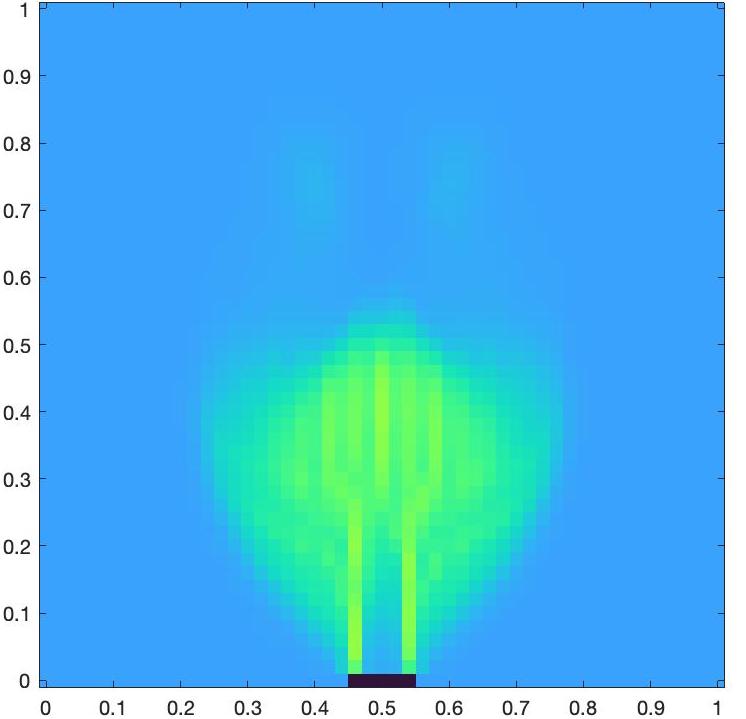}}
\hspace{0 mm}
\subfigure[][$\theta=0.25$]{\includegraphics[width=3.0cm, height=3.0cm]{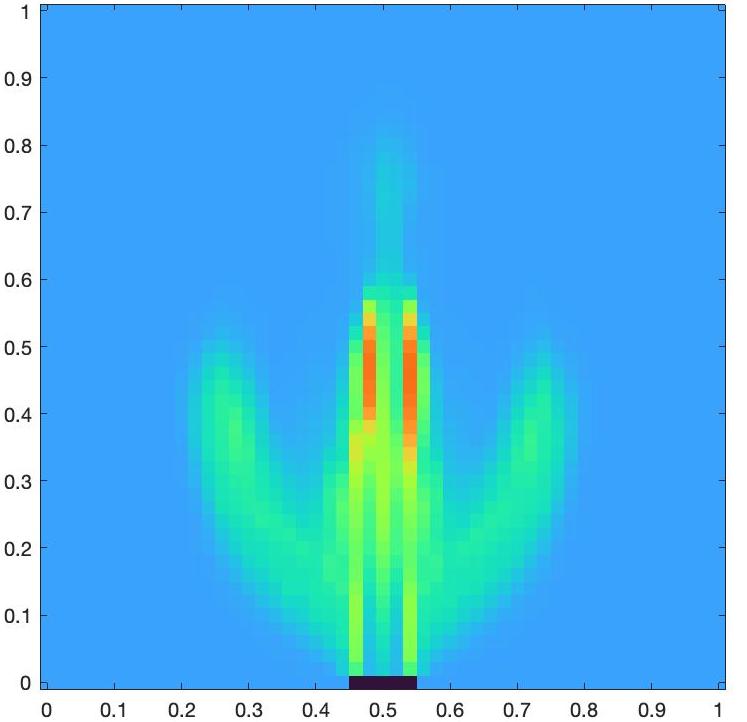}}
\hspace{0 mm}
\subfigure[][$\theta=2.5$]{\includegraphics[width=3.3cm, height=3.0cm]{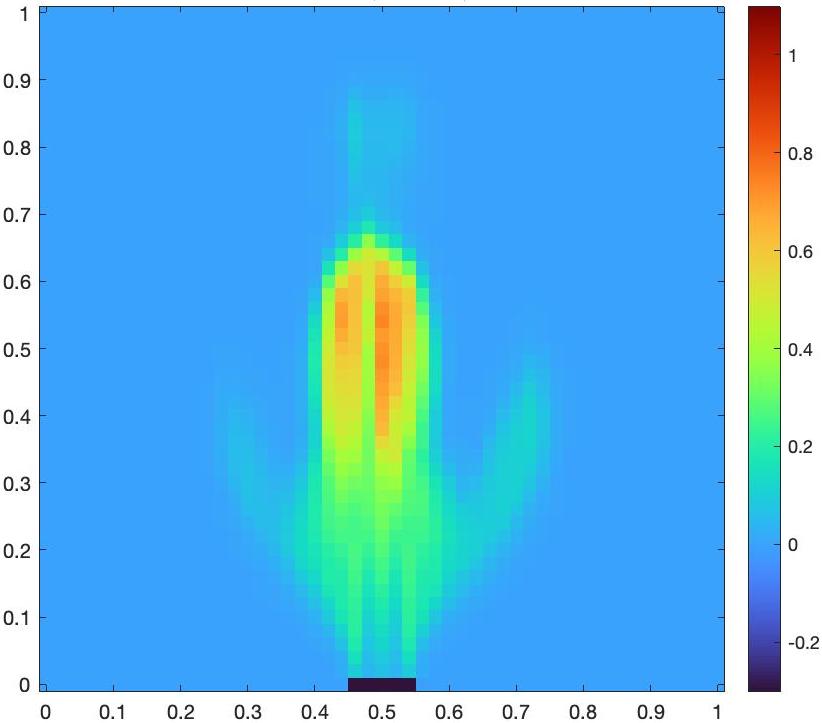}}
\caption{Test 5: Density evolution for three different values of $\theta$}
\label{fig:test5}
\end{figure}

Fig.\ \ref{fig:test5} shows three frames of the density evolution for three different values of $\theta$, respectively; the parameters $\bar t$ and $\Theta$ are the same for all cases and are set to $\bar t=0.48$ and $\Theta=0.24$, thus $\bar t-\Theta=0.24$.

(a) In the case $\theta=0$, the crowd initially moves upward toward the target. As usual in the Hughes's model, the part of the crowd furthest from the exit moves a bit to the side trying to overcome the people in front. This is indeed the best strategy if the rest of crowd stood still (which they do not). As a consequence, the typical heart-shaped configuration appears. 
(d) At $t=0.26>\bar t-\Theta$ the crowd is already aware that the exit will change position. Nevertheless, people are not able to forecast the formation of the congestion near the exit on top and then find it convenient not to change direction. After a while, though, they see the formation of the queue and turn down. 
(g) Now pedestrians which did not left the domain yet head down to reach the exit at the bottom.

In the case $\theta=T$ (MFG), the behaviour is rather different: 
(c) the heart-shaped configuration is less evident since people, being able forecasting the others, understand that moving aside is not convenient. 
Moreover, at time $\bar t-\Theta$ (when the crowd learns of the forthcoming change of the exit), people are already able to forecast the forthcoming congestion at the top, then the crowd immediately splits and a part of it starts moving downward. 
(f) At $t=0.26>\bar t-\Theta$ the crowd is already separated and a part of it is moving down.
(i) Finally, all the remaining pedestrians move down and leave the domain from the bottom.

(b,e,h) The case $\theta=0.25$ is in between the two extreme assumptions of the other models and seems to achieve a compromise regarding the behaviour of the crowd. 
We can see that the crowd both widens to the sides (as real crowds do) and splits as soon as people become aware of the forthcoming change.  


\section*{Funding}
This work was carried out within the research project ``SMARTOUR: Intelligent Platform for Tourism" (No. SCN\_00166) funded by the Ministry of University and Research with the Regional Development Fund of European Union (PON Research and Competitiveness 2007–2013). 

The authors also acknowledge the Italian Minister of Instruction, University and Research for supporting this research with funds coming from the project entitled Innovative numerical methods for evolutionary partial differential equations and applications (PRIN Project 2017, No. 2017KKJP4X). 

E.C.\ and M.M. are members of the INdAM Research group GNCS.

\section*{Acknowledgments}
The authors want to thank Fabio Camilli and Simone Cacace for the valuable help.

\end{document}